\input amstex.tex
\documentstyle{amsppt}
\magnification=\magstep1
\vsize 8.5truein
\hsize 6truein

\define\flow{\left(\bold{M},\{S^t\}_{t\in\Bbb R},\mu\right)}

\define\traj{S^{[a,b]}x_0}
\define\endrem{}
\define\harmas{\left(\Sigma,\Cal A,\vec\tau\right)}
\define\harmasv{\left(\Sigma',\Cal A',\vec\tau'\right)}

\noindent
November 27, 2003

\bigskip \bigskip

\heading
Proof of the Ergodic Hypothesis \\
for Typical Hard Ball Systems
\endheading
 
\bigskip \bigskip
 
\centerline{{\bf N\'andor Sim\'anyi}
\footnote{Research supported by the National Science Foundation, grant
DMS-0098773.}}

\bigskip \bigskip

\centerline{University of Alabama at Birmingham}
\centerline{Department of Mathematics}
\centerline{Campbell Hall, Birmingham, AL 35294 U.S.A.}
\centerline{E-mail: simanyi\@math.uab.edu}

\bigskip \bigskip

\hbox{\centerline{\vbox{\hsize 8cm {\bf Abstract.} We consider the system of
$N$ ($\ge2$) hard balls with masses $m_1,\dots,m_N$ and radius $r$ in the 
flat torus $\Bbb T_L^\nu=\Bbb R^\nu/L\cdot\Bbb Z^\nu$ of size $L$, $\nu\ge3$.
We prove the ergodicity (actually, the Bernoulli mixing property) of such 
systems for almost every selection $(m_1,\dots,m_N;\,L)$ of the outer 
geometric parameters. This theorem complements my earlier result that proved
the same, almost sure ergodicity for the case $\nu=2$. The method of that 
proof was primarily dynamical-geometric, whereas the present approach is
inherently algebraic.}}}

\bigskip \bigskip

\noindent
Primary subject classification: 37D50

\medskip

\noindent
Secondary subject classification: 34D05

\bigskip \bigskip

\heading
\S1. Introduction
\endheading

\bigskip \bigskip

Hard ball systems or, a bit more generally, mathematical billiards
constitute an important and quite interesting family of dynamical systems
being intensively studied by dynamicists and researchers of mathematical
physics, as well. These dynamical systems pose many challenging mathematical
questions, most of them concerning the ergodic (mixing) properties of such
systems. The introduction of hard ball systems and the first major steps in
their investigations date back to the 40's and 60's, see Krylov's paper
[K(1942)] and Sinai's ground-breaking works [Sin(1963)] and [Sin(1970)], 
in which the author --- among other things --- formulated the modern version
of Boltzmann's ergodic hypothesis (what we call today the Boltzmann--Sinai
ergodic hypothesis) by claiming that every hard ball system in a flat torus
is ergodic, of course after fixing the values of the trivial flow-invariant
quantities. In the articles
[Sin(1970)] and [B-S(1973)] Bunimovich and Sinai proved
this hypothesis for two hard disks on the two-dimensional unit torus
$\Bbb T^2$. The generalization of this result to higher dimensions $\nu>2$
took fourteen years, and was done by Chernov and Sinai in [S-Ch(1987)].
Although the model of two hard balls in $\Bbb T^\nu$ is already
rather involved technically, it is still a so called strictly dispersive
billiard system, i. e. such that the smooth components of the boundary
$\partial\bold Q$ of the configuration space are strictly concave from
outside $\bold Q$. (They are bending away from $\bold Q$.)
The billiard systems of more than two hard balls in
$\Bbb T^\nu$ are no longer strictly dispersive, but just semi-dispersive
(strict concavity of the smooth components of $\partial\bold Q$
is lost, merely concavity persists), and this circumstance causes a lot
of additional technical troubles in their study. In the series of my joint
papers with A. Kr\'amli and D. Sz\'asz [K-S-Sz(1989)], [K-S-Sz(1990)],
[K-S-Sz(1991)], and  [K-S-Sz(1992)]
we developed several new methods, and proved the ergodicity of 
more and more complicated semi-dispersive billiards culminating in the proof
of ergodicity of four billiard balls in the torus $\Bbb T^\nu$ 
($\nu\ge 3$), [K-S-Sz(1992)]. Then, in 1992, Bunimovich, Liverani,
Pellegrinotti and Sukhov [B-L-P-S(1992)] were able to prove the ergodicity for
some systems with an arbitrarily large number of hard balls. The shortcoming
of their model, however, is that, on one hand, they restrict the types of all
feasible ball-to-ball collisions, on the other hand they introduce some
additional
scattering effect with the collisions at the strictly concave wall of the
container. The only result with an arbitrarily large number of balls in a
flat unit torus $\Bbb T^\nu$ was achieved in the twin papers of mine
[Sim(1992-I-II)], where I managed to
prove the ergodicity (actually, the K-mixing property) of $N$ hard balls in
$\Bbb T^\nu$, provided that $N\le\nu$. The annoying shortcoming of that result
is that the larger the number of balls $N$ is, larger and larger dimension
$\nu$ of the ambient container is required by the method of the proof.

On the other hand, if someone considers a hard ball system in an
elongated torus which is long in one direction but narrow in the others,
so that the balls must keep their cyclic order in the ``long direction''
(Sinai's ``pen-case'' model), then the technical difficulties can be handled,
thanks to the fact that the collisions of balls are now restricted to
neighboring pairs. The hyperbolicity of such models in three dimensions and
the ergodicity in dimension four have been proved in [S-Sz(1995)]. 

The positivity of the metric entropy for several systems of hard balls
can be proven relatively easily, as was shown in the paper [W(1988)]. 
The articles [L-W(1995)] and [W(1990)] are nice surveys describing a 
general setup leading to the technical problems treated in a series of
research papers. For a comprehensive survey of the results and open problems 
in this field, see [Sz(1996)].

Pesin's theory [P(1977)] on the ergodic properties of non-uniformly hyperbolic,
smooth dynamical systems has been generalized substantially to dynamical 
systems with singularities (and with a relatively mild behavior near the
singularities) by A. Katok and J-M. Strelcyn [K-S(1986)]. Since then, the
so called Pesin's and Katok-Strelcyn's theories have become part of the 
folklore in the theory of dynamical systems. They claim that --- under some 
mild regularity conditions, particularly near the singularities --- every
non-uniformly hyperbolic and ergodic flow enjoys the Kolmogorov-mixing
property, shortly the K-mixing property.

Later on it was discovered and proven in [C-H(1996)] and [O-W(1998)] that the
above mentioned fully hyperbolic and ergodic flows with singularities turn out
to be automatically having the Bernoulli mixing (B-mixing) property. It is
worth noting here that almost every semi-dispersive billiard system, 
especially every hard ball system, enjoys those mild regularity conditions
imposed on the systems (as axioms) by [K-S(1986)], [C-H(1996)], and
[O-W(1998)]. In other words, for a hard ball flow 
$\left(\bold M,\{S^t\},\mu\right)$ the (global) ergodicity of the system
actually implies its full hyperbolicity and the B-mixing property, as well.

Finally, in our joint venture with D. Sz\'asz [S-Sz(1999)], we prevailed
over the difficulty caused by the low value of the dimension $\nu$ by
developing a brand new algebraic approach for the study of hard ball systems.
That result, however, only establishes complete hyperbolicity (nonzero Lyapunov
exponents almost everywhere) for $N$ balls in $\Bbb T^\nu$. The ergodicity
appeared to be a harder task.

We note, however, that the algebraic method developed in [S-Sz(1999)] 
is being further developed in this paper in order to obtain ergodicity,
not only full hyperbolicity.

Consider the $\nu$-dimensional ($\nu\ge2$), standard, flat torus
$\Bbb T_L^\nu=\Bbb R^\nu/L\cdot\Bbb Z^\nu$ as the vessel containing 
$N$ ($\ge2$) hard balls (spheres) $B_1,\dots,B_N$ with positive masses 
$m_1,\dots,m_N$ and (just for simplicity) common radius $r>0$. We always
assume that the radius $r>0$ is not too big, so
that even the interior of the arising
configuration space $\bold Q$ is connected. Denote the center of the ball
$B_i$ by $q_i\in\Bbb T^\nu$, and let $v_i=\dot q_i$ be the velocity of the
$i$-th particle. We investigate the uniform motion of the balls
$B_1,\dots,B_N$ inside the container $\Bbb T^\nu$ with half a unit of total 
kinetic energy: $E=\dfrac{1}{2}\sum_{i=1}^N m_i||v_i||^2=\dfrac{1}{2}$.
We assume that the collisions between balls are perfectly elastic. Since
--- beside the kinetic energy $E$ --- the total momentum
$I=\sum_{i=1}^N m_iv_i\in\Bbb R^\nu$ is also a trivial first integral of the
motion, we make the standard reduction $I=0$. Due to the apparent translation
invariance of the arising dynamical system, we factorize the configuration
space with respect to uniform spatial translations as follows:
$(q_1,\dots,q_N)\sim(q_1+a,\dots,q_N+a)$ for all translation vectors
$a\in\Bbb T^\nu$. The configuration space $\bold Q$ of the arising flow
is then the factor torus
$\left(\left(\Bbb T^\nu\right)^N/\sim\right)\cong\Bbb T^{\nu(N-1)}$
minus the cylinders

$$
C_{i,j}=\left\{(q_1,\dots,q_N)\in\Bbb T^{\nu(N-1)}\colon\;
\text{dist}(q_i,q_j)<2r \right\}
$$
($1\le i<j\le N$) corresponding to the forbidden overlap between the $i$-th
and $j$-th spheres. Then it is easy to see that the compound 
configuration point

$$
q=(q_1,\dots,q_N)\in\bold Q=\Bbb T^{\nu(N-1)}\setminus
\bigcup_{1\le i<j\le N}C_{i,j}
$$
moves in $\bold Q$ uniformly with unit speed and bounces back from the
boundaries $\partial C_{i,j}$ of the cylinders $C_{i,j}$ according to the
classical law of geometric optics: the angle of reflection equals the angle of
incidence. More precisely: the post-collision velocity $v^+$ can be obtained
from the pre-collision velocity $v^-$ by the orthogonal reflection across the
tangent hyperplane of the boundary $\partial\bold Q$ at the point of collision.
Here we must emphasize that the phrase ``orthogonal'' should be understood 
with respect to the natural Riemannian metric (the so called mass metric)
$||dq||^2=\sum_{i=1}^N m_i||dq_i||^2$ in the configuration space $\bold Q$.
For the normalized Liouville measure $\mu$ of the arising flow
$\{S^t\}$ we obviously have $d\mu=\text{const}\cdot dq\cdot dv$, where
$dq$ is the Riemannian volume in $\bold Q$ induced by the above metric and
$dv$ is the surface measure (determined by the restriction of the
Riemannian metric above) on the sphere of compound velocities

$$
\Bbb S^{\nu(N-1)-1}=\left\{(v_1,\dots,v_N)\in\left(\Bbb R^\nu\right)^N\colon\;
\sum_{i=1}^N m_iv_i=0 \text{ and } \sum_{i=1}^N m_i||v_i||^2=1 \right\}.
$$
The phase space $\bold M$ of the flow $\{S^t\}$ is the unit tangent bundle
$\bold Q\times\Bbb S^{d-1}$ of the configuration space $\bold Q$. (We will 
always use the shorthand notation $d=\nu(N-1)$ for the dimension of the 
billiard table $\bold Q$.) We must, however, note here that at the boundary
$\partial\bold Q$ of $\bold Q$ one has to glue together the pre-collision and
post-collision velocities in order to form the phase space $\bold M$, so
$\bold M$ is equal to the unit tangent bundle $\bold Q\times\Bbb S^{d-1}$
modulo this identification.

A bit more detailed definition of hard ball systems with arbitrary masses,
as well as their role in the family of cylindric billiards, can be found in
\S4 of [S-Sz(2000)] and in \S1 of [S-Sz(1999)]. We denote the
arising flow by $\flow$.

In the series of articles [K-S-Sz(1989)], [K-S-Sz(1991)], 
[K-S-Sz(1992)], [Sim(1992-I)], and [Sim(1992-II)]
the authors developed a powerful, three-step strategy for 
proving the (hyperbolic) ergodicity of hard ball systems. First of all,
all these proofs are inductions on the number $N$ of balls involved in the
problem. Secondly, the induction step itself consists of the following three
major steps:

\medskip

\subheading{Step I} To prove that every non-singular (i. e. smooth)
trajectory segment $\traj$ with a ``combinatorially rich'' (in a well
defined sense) symbolic collision sequence is automatically sufficient
(or, in other words, ``geometrically hyperbolic'', see below in \S2), 
provided that the phase point $x_0$ does not belong to a countable union $J$
of smooth sub-manifolds with codimension at least two. (Containing the 
exceptional phase points.)

The exceptional set $J$ featuring this result is negligible in our dynamical
considerations --- it is a so called slim set. For the basic properties of
slim sets, see \S2 below.

\medskip

\subheading{Step II} Assume the induction hypothesis, i. e. that all hard
ball systems with $N'$ balls ($2\le N'<N$) are (hyperbolic and) ergodic. 
Prove that then there exists a slim set $S\subset\bold M$ 
(see \S2) with the following property:
For every phase point $x_0\in\bold M\setminus S$ the entire trajectory
$S^{\Bbb R}x_0$ contains at most one singularity and its symbolic collision
sequence is combinatorially rich, just as required by the result of Step I.

\medskip

\subheading{Step III} By using again the induction hypothesis, prove that
almost every singular trajectory is sufficient in the time interval
$(t_0,+\infty)$, where $t_0$ is the time moment of the singular reflection.
(Here the phrase ``almost every'' refers to the volume defined by the induced
Riemannian metric on the singularity manifolds.)

We note here that the almost sure sufficiency of the singular trajectories
(featuring Step III) is an essential condition for the proof of the celebrated
Theorem on Local Ergodicity for algebraic semi-dispersive billiards 
proved by B\'alint--Chernov--Sz\'asz--T\'oth in [B-Ch-Sz-T (2002)].
Under this assumption the theorem of [B-Ch-Sz-T (2002)] states that in any 
algebraic
semi-dispersive billiard system (i. e. in a system such that the smooth
components of the boundary $\partial\bold Q$ are algebraic hypersurfaces)
a suitable, open neighborhood $U_0$ of any sufficient phase point 
$x_0\in\bold M$ (with at most one singularity on its trajectory) belongs
to a single ergodic component of the billiard flow $\flow$. 

In an inductive proof of ergodicity, steps I and II together ensure that
there exists an arc-wise connected set
$C\subset\bold M$ with full measure, such that every phase point $x_0\in C$
is sufficient with at most one singularity on its trajectory. Then the cited
Theorem on Local Ergodicity (now taking advantage of the result of Step III)
states that for every phase point $x_0\in C$ an open neighborhood $U_0$ of
$x_0$ belongs to one ergodic component of the flow. Finally, the connectedness
of the set $C$ and $\mu(\bold M\setminus C)=0$ easily imply that the flow
$\flow$ (now with $N$ balls) is indeed ergodic, and actually fully hyperbolic,
as well.

\medskip

The main result of this paper is the

\medskip

\subheading{Theorem} In the case $\nu\ge3$
for almost every selection $(m_1,\dots,m_N;\,L)$ of the outer geometric
parameters from the region $m_i>0$, $L>L_0(r,\,\nu)$ where the interior 
of the phase space is connected it is true that the billiard flow
$\left(\bold M_{\vec m,L},\{S^t\},\mu_{\vec m,L}\right)$ of the $N$-ball
system is ergodic and completely hyperbolic. Then, following from the results
of Chernov--Haskell [C-H(1996)] and Ornstein--Weiss [O-W(1998)], such a
semi-dispersive billiard system actually enjoys the B-mixing property, as well.

\medskip

\subheading{Remark 1} We note that the main result
of this paper and that of [Sim(2003)] nicely complement each other. They
precisely assert the same, almost sure ergodicity of hard ball systems in the
cases $\nu\ge3$ and $\nu=2$, respectively. It should be noted, however, that
the proof of [Sim(2003)] is primarily dynamical--geometric (except the 
verification of the Chernov-Sinai Ansatz), whereas the novel parts of the 
present proof are fundamentally algebraic.

\medskip

\subheading{Remark 2} The above inequality $L>L_0(r,\,\nu)$ corresponds to 
physically relevant situations. Indeed, in the case $L<L_0(r,\,\nu)$ the 
particles would not have enough room even to freely exchange positions. 

\medskip

The paper is organized as follows: \S2 provides all necessary prerequisites
and technical tools that will be required by the proof of the theorem. 
Based on the results obtained in [S-Sz(1999)],
the subsequent section \S3 carries out
Step I of the inductive strategy outlined above, but for the case when the
outer geometric parameters $(m_1,\dots,m_N;\,L)$ are incorporated in the
algebraic process as variables. (Just as the positions and velocities of 
the particles!) Finally, the closing section \S4 utilizes a ``Fubini type
argument'' by proving Step I for almost every (with respect to the Lebesgue
measure of the $(m_1,\dots,m_N;\,L)$-space) hard ball system ($N\ge2$,
$\nu\ge3$). This will finish the inductive proof of the Theorem, for Steps
II and III of the induction strategy are easy consequences of some earlier
results.

\bigskip \bigskip

\heading
2. Prerequisites
\endheading
 
\bigskip \bigskip
 
\subheading{2.1 Cylindric billiards} Consider the $d$-dimensional
($d\ge2$) flat torus $\Bbb T^d=\Bbb R^d/\Cal L$ supplied with the
usual Riemannian inner product $\langle\, .\, ,\, .\, \rangle$ inherited
from the standard inner product of the universal covering space $\Bbb R^d$.
Here $\Cal L\subset\Bbb R^d$ is assumed to be a lattice, i. e. a discrete
subgroup of the additive group $\Bbb R^d$ with $\text{rank}(\Cal L)=d$.
The reason why we want to allow general lattices, other than just the
integer lattice $\Bbb Z^d$, is that otherwise the hard ball systems would
not be covered. The geometry of the structure lattice $\Cal L$ in the
case of a hard ball system is significantly different from the geometry
of the standard lattice $\Bbb Z^d$ in the standard Euclidean space
$\Bbb R^d$, see later in this section.

The configuration space of a cylindric billiard is
$\bold Q=\Bbb T^d\setminus\left(C_1\cup\dots\cup C_k\right)$, where the
cylindric scatterers $C_i$ ($i=1,\dots,k$) are defined as follows.

Let $A_i\subset\Bbb R^d$ be a so called lattice subspace of $\Bbb R^d$,
which means that $\text{rank}(A_i\cap\Cal L)=\text{dim}A_i$. In this case
the factor $A_i/(A_i\cap\Cal L)$ is a sub-torus in $\Bbb T^d=\Bbb R^d/\Cal L$
which will be taken as the generator of the cylinder 
$C_i\subset\Bbb T^d$, $i=1,\dots,k$. Denote by $L_i=A_i^\perp$ the
ortho-complement of $A_i$ in $\Bbb R^d$. Throughout this paper we will
always assume that $\text{dim}L_i\ge2$. Let, furthermore, the numbers
$r_i>0$ (the radii of the spherical cylinders $C_i$) and some translation
vectors $t_i\in\Bbb T^d=\Bbb R^d/\Cal L$ be given. The translation
vectors $t_i$ play a role in positioning the cylinders $C_i$
in the ambient torus $\Bbb T^d$. Set
$$
C_i=\left\{x\in\Bbb T^d\colon\; \text{dist}\left(x-t_i,A_i/(A_i\cap\Cal L)
\right)<r_i \right\}.
$$
In order to avoid further unnecessary complications, we always assume that
the interior of the configuration space 
$\bold Q=\Bbb T^d\setminus\left(C_1\cup\dots\cup C_k\right)$ is connected.
The phase space $\bold M$ of our cylindric billiard flow will be the
unit tangent bundle of $\bold Q$ (modulo the natural gluing at its
boundary), i. e. $\bold M=\bold Q\times\Bbb S^{d-1}$. (Here $\Bbb S^{d-1}$
denotes the unit sphere of $\Bbb R^d$.)

The dynamical system $\flow$, where $S^t$ ($t\in\Bbb R$) is the dynamics 
defined by the uniform motion inside the domain $\bold Q$ and specular
reflections at its boundary (at the scatterers), and $\mu$ is the
Liouville measure, is called a cylindric billiard flow.

We note that the cylindric billiards --- defined above --- belong to the wider
class of so called semi-dispersive billiards, which means that the smooth 
components $\partial\bold Q_i$ of the boundary $\partial\bold Q$
of the configuration space $\bold Q$ are (not necessarily strictly) concave
from outside of $\bold Q$,
i. e. they are bending away from the interior of $\bold Q$. As to the notions 
and notations in connection with semi-dispersive billiards, the reader is 
kindly referred to the article [K-S-Sz(1990)].

\medskip

Throughout this paper we will always assume --- without explicitly stating
--- that the considered semi-dispersive billiard system fulfills the
following conditions:

\medskip

$$
\text{int}\bold Q\text{ is connected, and}
\tag 2.1.1
$$

$$
\aligned
&\text{the }d\text{-dim spatial angle }\alpha(q)\text{ subtended by }\bold Q \\
&\text{at any of its boundary points }q\in\partial\bold Q\text{ is uniformly
positive.}
\endaligned
\tag 2.1.2
$$
We note, however, that in the case of hard ball systems with a fixed radius
$r$ of the balls (see below) the
non-degeneracy condition (2.1.2) only excludes countably many values of the
size $L$ of the container torus $\Bbb T^\nu_L=\Bbb R^\nu/L\cdot\Bbb Z^\nu$
from the region $L>L_0(r,\,\nu)$ where (2.1.1) is true. Therefore, in the
sense of our theorem of ``almost sure ergodicity'', the non-degeneracy condition
(2.1.2) does not mean a restriction of generality.

\bigskip

\subheading{2.2 Hard ball systems} Hard ball systems in the flat
torus $\Bbb T^\nu_L=\Bbb R^\nu/L\cdot\Bbb Z^\nu$ ($\nu\ge2$) with positive masses
$m_1,\dots,m_N$ are described (for example) in \S 1 of [S-Sz(1999)].
These are the dynamical systems describing the motion of $N$ ($\ge2$) hard
balls with a common radius $r>0$ and positive masses $m_1,\dots,m_N$ in 
the flat torus of size $L$, $\Bbb T^\nu_L=\Bbb R^\nu/L\cdot\Bbb Z^\nu$.
(Just for simplicity, we will assume that the radii have the common 
value $r$.) The center of the
$i$-th ball is denoted by $q_i$ ($\in\Bbb T^\nu_L$), its time derivative is
$v_i=\dot q_i$, $i=1,\dots,N$. One uses the standard reduction of kinetic
energy $E=\frac{1}{2}\sum_{i=1}^N m_i||v_i||^2=\frac{1}{2}$.
The arising configuration space (still without the removal of the scattering
cylinders $C_{i,j}$) is the torus

$$
\Bbb T_L^{\nu N}=\left(\Bbb T_L^{\nu}\right)^N=\left\{(q_1,\dots,q_N)\colon\;
q_i\in\Bbb T_L^\nu,\; i=1,\dots,N\right\}
$$
supplied with the Riemannian inner product (the so called mass metric)

$$
\langle v,v'\rangle=\sum_{i=1}^N m_i\langle v_i,v'_i \rangle
\tag 2.2.1
$$
in its common tangent space $\Bbb R^{\nu N}=\left(\Bbb R^{\nu}\right)^N$.
Now the Euclidean space $\Bbb R^{\nu N}$ with the inner product (2.2.1)
plays the role of $\Bbb R^d$ in the original definition of cylindric
billiards, see \S 2.1 above.

The generator subspace $A_{i,j}\subset \Bbb R^{\nu N}$ ($1\le i<j\le N$)
of the cylinder $C_{i,j}$ (describing the collisions between the $i$-th and
$j$-th balls) is given by the equation
$$
A_{i,j}=\left\{(q_1,\dots,q_N)\in\left(\Bbb R^\nu\right)^N\colon\;
q_i=q_j \right\},
\tag 2.2.2
$$
see (4.3) in [S-Sz(2000)]. Its ortho-complement 
$L_{i,j}\subset\Bbb R^{\nu N}$ is then defined by the equation
$$
L_{i,j}=\left\{(q_1,\dots,q_N)\in\left(\Bbb R^\nu\right)^N\colon\;
q_k=0 \text{ for } k\ne i,j, \text{ and } m_iq_i+m_jq_j=0 \right\},
\tag 2.2.3
$$
see (4.4) in [S-Sz(2000)].
Easy calculation shows that the cylinder $C_{i,j}$ 
(describing the overlap of the $i$-th and $j$-th balls)
is indeed spherical and the radius of its base sphere is equal to
$r_{i,j}=2r\sqrt{\frac{m_im_j}{m_i+m_j}}$, see \S 4, especially formula
(4.6) in [S-Sz(2000)].

The structure lattice $\Cal L\subset\Bbb R^{\nu N}$ is clearly the
lattice $\Cal L=\left(L\cdot\Bbb Z^{\nu}\right)^N=L\cdot\Bbb Z^{N\nu}$. 

Due to the presence of an extra invariant quantity
$I=\sum_{i=1}^N m_iv_i$, one usually makes the reduction
$\sum_{i=1}^N m_iv_i=0$ and, correspondingly, factorizes the configuration
space with respect to uniform spatial translations:

$$
(q_1,\dots,q_N)\sim(q_1+a,\dots,q_N+a), \quad a\in\Bbb T_L^\nu.
\tag 2.2.4
$$
The natural, common tangent space of this reduced configuration space is then

$$
\Cal Z=\left\{(v_1,\dots,v_N)\in\left(\Bbb R^\nu\right)^N\colon\;
\sum_{i=1}^N m_iv_i=0\right\}=\left(\bigcap_{i<j}A_{i,j}
\right)^\perp=\left(\Cal A\right)^\perp
\tag 2.2.5
$$
supplied again with the inner product (2.2.1), see also (4.1) and
(4.2) in [S-Sz(2000)]. The base spaces $L_{i,j}$
of (2.2.3) are obviously subspaces of $\Cal Z$, and we take
$\tilde A_{i,j}=A_{i,j}\cap\Cal Z=P_{\Cal Z}(A_{i,j})$ as the ortho-complement
of $L_{i,j}$ in $\Cal Z$. (Here $P_{\Cal Z}$ denotes the orthogonal projection
onto the space $\Cal Z$.)

Note that the configuration space of the reduced system (with 
the identification (2.2.4)) is naturally the torus
$\Bbb R^{\nu N}/(\Cal A+L\cdot\Bbb Z^{\nu N})=
\Cal Z/P_{\Cal Z}(L\cdot\Bbb Z^{\nu N})$.

\bigskip

\subheading{2.3. Collision graphs} Let $S^{[a,b]}x$ be a nonsingular, finite
trajectory segment with the collisions $\sigma_1,\dots,\sigma_n$ listed in
time order. (Each $\sigma_k$ is an unordered pair $(i,j)$ of different labels
$i,j\in\{1,2,\dots,N\}$.) The graph $\Cal G=(\Cal V,\Cal E)$ with vertex set
$\Cal V=\{1,2,\dots,N\}$ and set of edges
$\Cal E=\{\sigma_1,\dots,\sigma_n\}$ is called the {\it collision graph}
of the orbit segment $S^{[a,b]}x$. For a given positive number $C$, the
collision graph $\Cal G=(\Cal V,\Cal E)$ of the orbit segment $S^{[a,b]}x$
will be called {\it $C$-rich} if $\Cal G$ contains at least $C$ connected,
consecutive (i. e. following one after the other in time, according to the
time-ordering given by the trajectory segment $S^{[a,b]}x$) subgraphs. 

\bigskip

\subheading{2.4 Trajectory Branches}
We are going to briefly describe the discontinuity of the flow
$\{S^t\}$ caused by a multiple collisions at time $t_0$.
Assume first that the pre--collision velocities of the particles are given.
What can we say about the possible post--collision velocities? Let us perturb
the pre--collision phase point (at time $t_0-0$) infinitesimally, so that the
collisions at $\sim t_0$ occur at infinitesimally different moments. By
applying the collision laws to the arising finite sequence of collisions, we
see that the post-collision velocities are fully determined by the time--
ordering of the considered collisions. Therefore, the collection of all
possible time-orderings of these collisions gives rise to a finite family of
continuations of the trajectory beyond $t_0$. They are called the
trajectory branches. It is quite clear that similar statements can be
said regarding the evolution of a trajectory through a multiple collision
in reverse time. Furthermore, it is also obvious that for any given
phase point $x_0\in\bold M$ there are two, $\omega$-high trees
$\Cal T_+$ and $\Cal T_-$ such that $\Cal T_+$ ($\Cal T_-$) describes all the
possible continuations of the positive (negative) trajectory
$S^{[0,\infty)}x_0$ ($S^{(-\infty,0]}x_0$). (For the definitions of trees and
for some of their applications to billiards, cf. the beginning of \S 5
in [K-S-Sz(1992)].) It is also clear that all possible continuations
(branches) of the whole trajectory $S^{(-\infty,\infty)}x_0$ can be uniquely
described by all pairs $(B_-,B_+)$ of $\omega$-high branches of
the trees $\Cal T_-$ and $\Cal T_+$ ($B_-\subset\Cal T_-, B_+\subset
\Cal T_+$).

Finally, we note that the trajectory of the phase point $x_0$ has exactly two
branches, provided that $S^tx_0$ hits a singularity for a single value
$t=t_0$, and the phase point $S^{t_0}x_0$ does not lie on the intersection
of more than one singularity manifolds. In this case we say that the 
trajectory of $x_0$ has a ``simple singularity''.

\bigskip

\subheading{2.5 Neutral Subspaces, Advance, and Sufficiency}
Consider a nonsingular trajectory segment $S^{[a,b]}x$.
Suppose that $a$ and $b$ are not moments of collision.

\medskip

\proclaim{Definition 2.5.1} The neutral space $\Cal N_0(S^{[a,b]}x)$
of the trajectory segment $S^{[a,b]}x$ at time zero ($a<0<b$)\ is
defined by the following formula:

$$
\aligned
&\Cal N_0(S^{[a,b]}x)=\big \{W\in\Cal Z\colon\;\exists (\delta>0) \;
\text{ such that } \; \forall \alpha \in (-\delta,\delta) \\
&V\left(S^a\left(Q(x)+\alpha W,V(x)\right)\right)=V(S^ax)\text{ and }
V\left(S^b\left(Q(x)+\alpha W,V(x)\right)\right)=V(S^bx)\big\}.
\endaligned
$$
\endproclaim \endrem

\noindent
($\Cal Z$ is the common tangent space $\Cal T_q\bold Q$ of the parallelizable
manifold $\bold Q$ at any of its points $q$, while $V(x)$ is the velocity
component of the phase point $x=\left(Q(x),\,V(x)\right)$.)

It is known (see (3) in \S 3 of [S-Ch (1987)]) that
$\Cal N_0(S^{[a,b]}x)$ is a linear subspace of $\Cal Z$ indeed, and
$V(x)\in \Cal N_0(S^{[a,b]}x)$. The neutral space $\Cal N_t(S^{[a,b]}x)$
of the segment $S^{[a,b]}x$ at time $t\in [a,b]$ is defined as follows:

$$
\Cal N_t(S^{[a,b]}x)=\Cal N_0\left(S^{[a-t,b-t]}(S^tx)\right).
$$
It is clear that the neutral space $\Cal N_t(S^{[a,b]}x)$ can be
canonically
identified with $\Cal N_0(S^{[a,b]}x)$ by the usual identification of the
tangent spaces of $\bold Q$ along the trajectory $S^{(-\infty,\infty)}x$
(see, for instance, \S 2 of [K-S-Sz(1990)]).

Our next  definition is  that of  the advance. Consider a
non-singular orbit segment $S^{[a,b]}x$ with the symbolic collision sequence
$\Sigma=(\sigma_1, \dots, \sigma_n)$ ($n\ge 1$), meaning that $S^{[a,b]}x$
has exactly $n$ collisions with $\partial\bold Q$, and the $i$-th collision
($1\le i\le n$) takes place at the boundary of the cylinder $C_{\sigma_i}$.
For $x=(Q,V)\in\bold M$ and $W\in\Cal Z$, $\Vert W\Vert$ sufficiently small, 
denote $T_W(Q,V):=(Q+W,V)$.

\proclaim{Definition 2.5.2}
For any $1\le k\le n$ and $t\in[a,b]$, the advance
$$
\alpha(\sigma_k)\colon\;\Cal N_t(S^{[a,b]}x) \rightarrow  \Bbb R
$$
of the collision $\sigma_k$ is the unique linear extension of the linear 
functional $\alpha(\sigma_k)$
defined in a sufficiently small neighborhood of the origin of 
$\Cal N_t(S^{[a,b]}x)$ in the following way:
$$
\alpha(\sigma_k)(W):= t_k(x)-t_k(S^{-t}T_WS^tx).
$$
\endproclaim \endrem
Here $t_k=t_k(x)$ is the time moment of the $k$-th collision $\sigma_k$ on
the trajectory of $x$ after time $t=a$. The above formula and the notion of
the advance functional 

$$
\alpha_k=\alpha(\sigma_k):\; \Cal N_t\left(S^{[a,b]}x\right)\to\Bbb R
$$
has two important features:

\medskip

(i) If the spatial translation $(Q,V)\mapsto(Q+W,V)$ is carried out at time
$t$, then $t_k$ changes linearly in $W$, and it takes place just 
$\alpha_k(W)$ units of time earlier. (This is why it is called ``advance''.)

\medskip

(ii) If the considered reference time $t$ is somewhere between $t_{k-1}$
and $t_k$, then the neutrality of $W$ with respect to $\sigma_k$ precisely 
means that

$$
W-\alpha_k(W)\cdot V(x)\in A_{\sigma_k},
$$
i. e. a neutral (with respect to the collision $\sigma_k$) spatial translation
$W$ with the advance $\alpha_k(W)=0$ means that the vector $W$ belongs to the
generator space $A_{\sigma_k}$ of the cylinder $C_{\sigma_k}$.

It is now time to bring up the basic notion of sufficiency 
(or, sometimes it is also called geometric hyperbolicity) of a
trajectory (segment). This is the utmost important necessary condition for
the proof of the fundamental theorem for algebraic semi-dispersive billiards, 
see Theorem 4.4 in [B-Ch-Sz-T(2002)].

\medskip

\proclaim{Definition 2.5.3}
\roster
\item
The nonsingular trajectory segment $S^{[a,b]}x$ ($a$ and $b$ are supposed not
to be moments of collision) is said to be sufficient if and only if
the dimension of $\Cal N_t(S^{[a,b]}x)$ ($t\in [a,b]$) is minimal, i.e.
$\text{dim}\ \Cal N_t(S^{[a,b]}x)=1$.
\item
The trajectory segment $S^{[a,b]}x$ containing exactly one singularity (a so 
called ``simple singularity'', see 2.4 above) is said to be sufficient if 
and only if both branches of this trajectory segment are sufficient.
\endroster
\endproclaim \endrem

\medskip

\proclaim{Definition 2.5.4}
The phase point $x\in\bold M$ with at most one (simple) singularity is said
to be sufficient if and only if its whole trajectory $S^{(-\infty,\infty)}x$
is sufficient, which means, by definition, that some of its bounded
segments $S^{[a,b]}x$ are sufficient.
\endproclaim \endrem

In the case of an orbit $S^{(-\infty,\infty)}x$ with a simple
singularity, sufficiency means that both branches of
$S^{(-\infty,\infty)}x$ are sufficient.

\bigskip

\subheading{2.6. No accumulation (of collisions) in finite time} 
By the results of Vaserstein [V(1979)], Galperin [G(1981)] and
Burago-Ferleger-Kononenko [B-F-K(1998)], in a semi-dis\-per\-sive billiard 
flow with the property (2.1.2) there can only be finitely many 
collisions in finite time intervals, see Theorem 1 in [B-F-K(1998)]. 
Thus, the dynamics is well defined as long as the trajectory does not hit 
more than one boundary components at the same time.

\bigskip

\subheading{2.7. Slim sets} 
We are going to summarize the basic properties of codimension-two subsets $A$
of a connected, smooth manifold $M$ with a possible boundary. Since these subsets 
$A$ are just those negligible in our dynamical discussions, we shall call them 
slim. As to a  broader exposition of the issues, see [E(1978)] or \S2 of
[K-S-Sz(1991)].

Note that the dimension $\dim A$ of a separable metric space $A$ is one of the
three classical notions of topological dimension: the covering 
(\v Cech-Lebesgue), the small inductive (Menger-Urysohn), or the large 
inductive (Brouwer-\v Cech) dimension. As it is known from general general 
topology, all of them are the same for separable metric spaces.

\medskip

\subheading{Definition 2.7.1}
A subset $A$ of $M$ is called slim if and only if $A$ can be covered by a 
countable family of codimension-two (i. e. at least two) closed sets of
$\mu$--measure zero, where $\mu$ is a smooth measure on $M$. (Cf.
Definition 2.12 of [K-S-Sz(1991)].)

\medskip

\subheading{Property 2.7.2} The  collection of all slim subsets of $M$ is a
$\sigma$-ideal, that is, countable unions of slim sets and arbitrary
subsets of slim sets are also slim.

\medskip

\subheading{Proposition 2.7.3. (Locality)}
A subset $A\subset M$ is slim if and only if for
every $x\in A$ there exists an open neighborhood $U$ of $x$ in $M$ such that 
$U\cap A$ is slim. (Cf. Lemma 2.14 of [K-S-Sz(1991)].)

\medskip

\subheading{Property 2.7.4} A closed subset $A\subset M$ is slim if and only
if $\mu(A)=0$ and $\dim A\le\dim M-2$.

\medskip

\subheading{Property 2.7.5. (Integrability)}
If $A\subset M_1\times M_2$ is a closed subset of the product of two smooth 
manifolds with possible boundaries, and for every $x\in M_1$ the set
$$
A_x=\{ y\in M_2\colon\; (x,y)\in A\}
$$
is slim in $M_2$, then $A$ is slim in $M_1\times M_2$.

\medskip

The following propositions characterize the codimension-one and 
codimension-two sets.

\medskip

\subheading{Proposition 2.7.6}
For any closed subset $S\subset M$ the following three conditions are 
equivalent:

\roster

\item"{(i)}" $\dim S\le\dim M-2$;

\item"{(ii)}"  $\text{int}S=\emptyset$ and for every open connected set 
$G\subset M$ the difference set $G\setminus S$ is also connected;

\item"{(iii)}" $\text{int}S=\emptyset$ and for every point $x\in M$ and for any
open neighborhood $V$ of $x$ in $M$ there exists a smaller open neighborhood
$W\subset V$ of the point $x$ such that for every pair of points 
$y,z\in W\setminus S$ there is a continuous curve $\gamma$ in the set 
$V\setminus S$ connecting the points $y$ and $z$.

\endroster

\noindent
(See Theorem 1.8.13 and Problem 1.8.E of [E(1978)].)

\medskip

\subheading{Proposition 2.7.7} For any subset $S\subset M$ the condition 
$\dim S\le\dim M-1$ is equivalent to $\text{int}S=\emptyset$.
(See Theorem 1.8.10 of [E(1978)].)

\medskip

We recall an elementary, but important lemma (Lemma 4.15 of [K-S-Sz(1991)]).
Let $R_2$ be the set of phase points 
$x\in\bold M\setminus\partial\bold M$ such that the trajectory 
$S^{(-\infty,\infty)}x$ has more than one singularities.

\subheading{Proposition 2.7.8}
The set $R_2$ is a countable union of codimension-two
smooth sub-manifolds of $M$ and, being such, it is slim.

\medskip

The next lemma establishes the most important property of slim sets which
gives us the fundamental geometric tool to connect the open ergodic components
of billiard flows.

\medskip

\subheading{Proposition 2.7.9}
If $M$ is connected, then the complement $M\setminus A$ of a slim $F_\sigma$ 
set $A\subset M$ is an arc-wise connected ($G_\delta$) set of full measure. 
(See Property 3 of \S 4.1 in [K-S-Sz(1989)]. The  $F_\sigma$ sets are, 
by definition, the countable unions of closed sets, while the $G_\delta$ sets
are the countable intersections of open sets.)

\medskip

\subheading{2.8. The subsets $\bold M^0$ and $\bold M^\#$} Denote by
$\bold M^\#$ the set of all phase points $x\in\bold M$ for which the
trajectory of $x$ encounters infinitely many non-tangential collisions
in both time directions. The trajectories of the points 
$x\in\bold M\setminus\bold M^\#$ are lines: the motion is linear and uniform,
see the appendix of [Sz(1994)]. It is proven in lemmas A.2.1 and A.2.2
of [Sz(1994)] that the closed set $\bold M\setminus\bold M^\#$ is a finite
union of hyperplanes. It is also proven in [Sz(1994)] that, locally, the two
sides of a hyper-planar component of $\bold M\setminus\bold M^\#$ can be
connected by a positively measured beam of trajectories, hence, from the point
of view of ergodicity, in this paper it is enough to show that the connected
components of $\bold M^\#$ entirely belong to one ergodic component. This is
what we are going to do in this paper.

Denote by $\bold M^0$ the set of all phase points $x\in\bold M^\#$ the 
trajectory of which does not hit any singularity, and use the notation
$\bold M^1$ for the set of all phase points $x\in\bold M^\#$ whose orbit
contains exactly one, simple singularity. According to Proposition 2.7.8,
the set $\bold M^\#\setminus(\bold M^0\cup\bold M^1)$ is a countable union of
smooth, codimension-two ($\ge2$) submanifolds of $\bold M$, and, therefore,
this set may be discarded in our study of ergodicity, please see also the
properties of slim sets above. Thus, we will restrict our attention to the
phase points $x\in\bold M^0\cup\bold M^1$.

\medskip

\subheading{2.9. The ``Chernov-Sinai Ansatz''} An essential precondition for
the Theorem on Local Ergodicity by B\'alint--Chernov--Sz\'asz--T\'oth 
(Theorem 4.4 of [B-Ch-Sz-T(2002)]) is the
so called ``Chernov-Sinai Ansatz'' which we are going to formulate below.
Denote by $\Cal S\Cal R^+\subset\partial\bold M$ the set of all phase points
$x_0=(q_0,v_0)\in\partial\bold M$ corresponding to singular reflections
(a tangential or a double collision at time zero) supplied with the 
post-collision (outgoing) velocity $v_0$. It is well known that
$\Cal S\Cal R^+$ is a compact cell complex with dimension
$2d-3=\text{dim}\bold M-2$. It is also known (see Lemma 4.1 in [K-S-Sz(1990)])
that for $\nu$-almost every phase point $x_0\in\Cal S\Cal R^+$ the forward 
orbit $S^{(0,\infty)}x_0$ does not hit any further singularity. 
(Here $\nu$ is the Riemannian volume of $\Cal S\Cal R^+$ induced by the 
restriction of the natural Riemannian metric of $\bold M$.)
The Chernov-Sinai Ansatz postulates that for $\nu$-almost every 
$x_0\in\Cal S\Cal R^+$ the forward orbit $S^{(0,\infty)}x_0$ is sufficient 
(geometrically hyperbolic). 

\medskip

\subheading{2.10. The Theorem on Local Ergodicity} The Theorem on Local 
Ergodicity by B\'alint--Chernov--Sz\'asz--T\'oth (Theorem 4.4 of
[B-Ch-Sz-T(2002)]) claims the following: Let $\flow$ be a semi-dispersive
billiard flow with (2.1.1)--(2.1.2) and with the property
that the smooth components of the boundary $\partial\bold Q$ of the 
configuration space are algebraic hyper-surfaces. (The cylindric billiards
automatically fulfill this algebraicity condition.) Assume -- further --
that the Chernov-Sinai Ansatz holds true, and a phase point 
$x_0\in\left(\bold M\setminus\partial\bold M\right)\cap\bold M^\#$ is 
given with the properties

\medskip

(i) $S^{(-\infty,\infty)}x$ has at most one singularity,

\noindent
and

(ii) $S^{(-\infty,\infty)}x$ is sufficient. 

\medskip

Then some open neighborhood $U_0\subset\bold M$ of $x_0$ belongs to a single
ergodic component of the flow $\flow$. (Modulo the zero sets, of course.)

\bigskip \bigskip

\heading
3. Non-sufficiency Occurs on a Codimension-two Set \\
The Case $\nu\ge3$
\endheading

\bigskip \bigskip

The opening part of this section contains a slightly modified version of
Lemma 4.43 from [S-Sz(1999)]. The reason why we had to modify the recursion
for the sequence $C(N)$ (from $C(N)=(N/2)\cdot\max\left\{C(N-1),\,3\right\}$ 
to $C(N)=(N/2)\cdot\left(2C(N-1)+1\right)$) is that our Corollary 3.5
(below) requires $(2C(N)+1)$-richness instead of the usual $C(N)$-richness.
In the present paper the sequence $C(N)$ always denotes the one defined by
the recursion in Lemma 3.1 instead of the one defined in Lemma 4.43 of 
[S-Sz(1999)]. This should not cause any confusion.

We note that the upcoming lemma is purely combinatorial.

\medskip

\subheading{Lemma 3.1}
Define the sequence of positive numbers $C(N)$
recursively by taking $C(2)=1$ and
$C(N)=(N/2)\cdot\left(2C(N-1)+1\right)$ for $N\ge3$.
Let $N\ge3$, and suppose that the symbolic collision sequence
$\Sigma=(\sigma_1,\dots,\sigma_n)$ for $N$ particles is $C(N)$-rich.
Then we can find a particle, say the one with label $N$, and two indices
$1\le p<q\le n$ such that

(i) $N\in\sigma_p\cap\sigma_q$,

(ii) $N\notin\bigcup_{j=p+1}^{q-1}\sigma_j$,

(iii) $\sigma_p=\sigma_q\Longrightarrow\,\left(\exists j\right)\;
\left(p<j<q\,\&\,\sigma_p\cap\sigma_j\ne\emptyset\right)$, and

(iv) $\Sigma'$ is $(2C(N-1)+1)$-rich on the vertex set $\{1,\dots,N-1\}$.

(Here, just as in the case of derived schemes, we denote by $\Sigma'$ the
symbolic sequence that can be obtained from $\Sigma$ by discarding all
edges containing $N$.)

\medskip

\subheading{Proof} The hypothesis on $\Sigma$ implies that there exist
subsequences $\Sigma_1,\dots,\Sigma_r$ of $\Sigma$ with the following
properties:

(1) For $1\le i<j\le r$ every collision of $\Sigma_i$ precedes every collision
of $\Sigma_j$,

(2) the graph of $\Sigma_i$ ($1\le i\le r$) is a tree (a connected graph
without loop) on the vertex set $\{1,\dots,N\}$, and

(3) $r\ge C(N)$.

\medskip

Since every tree contains at least two vertices with degree one and
$C(N)=(N/2)\cdot\left\{2C(N-1)+1\right\}$, there
is a vertex, say the one labeled by $N$, such that $N$ is a degree-one vertex
of $\Sigma_{i(1)},\dots,\Sigma_{i(t)}$, where $1\le i(1)<\dots<i(t)\le r$ and
$t\ge 2C(N-1)+1$. Thus (iv) obviously holds.

Let $\sigma_{p'}$ the edge of $\Sigma_{i(1)}$ that contains $N$ and, similarly,
let $\sigma_{q'}$ be the edge of $\Sigma_{i(t)}$ containing the vertex $N$.
Then the fact $t\ge 3$ ensures that the following properties hold:

(i)' $N\in\sigma_{p'}\cap\sigma_{q'}$,

(iii)' $\sigma_{p'}=\sigma_{q'}\Longrightarrow\, \exists j\; p'<j<q' \,\&\,
\sigma_{p'}\cap\sigma_j\ne\emptyset$, $\sigma_j\ne\sigma_{p'}$.

Let $\sigma_p$, $\sigma_q$ ($1\le p<q\le n$) be a pair of edges
$\sigma_{p'}$, $\sigma_{q'}$ ($1\le p'<q'\le n$) fulfilling (i)' and (iii)'
and having the minimum possible value of $q'-p'$. Elementary inspection shows
that then (ii) must also hold for $\sigma_p$, $\sigma_q$. Lemma 3.51 is now
proved. \qed

\medskip

Let us fix a triplet $\harmas$ of the discrete (combinatorial) orbit structure
with Property (A) (just as in [S-Sz(1999)], see Definition 3.31 there), 
and assume that
$\Sigma=(\sigma_1,\dots,\sigma_n)$ is $C(N)$-rich, i. e. it contains at least
$C(N)$ consecutive, connected collision graphs. We also consider the complex
analytic manifold $\Omega\harmas$ of all complex $\harmas$-orbits
$\omega$ (Definition 3.20 in [S-Sz(1999)]) and the open, dense, connected 
domain $D\harmas\subset\Bbb C^{(2\nu+1)N+1}$ of all allowable initial data
$\vec x=\vec x(\omega)$, see Definition 3.18 in [S-Sz(1999)].
Let, finally, $Q(\vec x)$ be a
common {\it irreducible} divisor of the polynomials
$P_1(\vec x),\dots,P_s(\vec x)$ from (4.3) in [S-Sz(1999)].
(If such a common divisor exists.)
In this section we will need several results about such common irreducible 
divisors $Q(\vec x)$ of the polynomials $P_1(\vec x),\dots,P_s(\vec x)$.

The first of them, as it is classically known from algebraic geometry (see,
for example, [M(1976)]), is that the solution set $V=\{Q(\vec x)=0\}$ of the 
equation $Q(\vec x)=0$ is a so called irreducible (or, indecomposable) complex
algebraic variety of codimension $1$ in $\Bbb C^{(2\nu+1)N+1}$, which means 
that $V$ is not the union of two, proper algebraic sub-varieties. Secondly, the 
smooth part $S$ of $V$ turns out to be a connected complex analytic manifold,
while the non-smooth part $V\setminus S$ of $V$ is a complex algebraic variety
of dimension strictly less that $\text{dim}V=(2\nu+1)N$, see [M(1976)].
Finally, if a polynomial $R(\vec x)$ vanishes on $V$, then $Q(\vec x)$ must
be a divisor of $R(\vec x)$. (The last statement is a direct consequence of
Hilbert's Theorem on Zeroes, see again [M(1976)].) 

The first result, specific to our current dynamical situation, is

\medskip

\subheading{Proposition 3.2} The polynomials $P_1(\vec x),\dots,P_s(\vec x)$
of (4.3) in [S-Sz(1999)] are homogeneous in the masses $m_1,\dots,m_N$ and, 
consequently, any common divisor $Q(\vec x)$ of these polynomials is also 
homogeneous in the masses.

\medskip

\subheading{Proof} It is clear that the complex dynamics encoded in the
orbits $\omega\in\Omega\harmas$ only depends on the ratios of masses
$m_2/m_1,\,m_3/m_1,\,\dots,\,m_N/m_1$. Consequently, all algebraic functions
$f_i(\vec x)$ ($i=1,\dots,s$) featuring the proof of Lemma 4.2 of
[S-Sz(1999)] are
homogeneous of degree $0$ in the masses. Since the rational function

$$
\frac{P_i(\vec x)}{Q_i(\vec x)}\in\Bbb K_0=\Bbb C(\vec x)
$$
is the product $\alpha$ of all conjugates of $f_i(\vec x)$ 
(see the proof of the lemma just cited), we get that 
$\dfrac{P_i(\vec x)}{Q_i(\vec x)}$ is also homogeneous of degree $0$ in the
variables $m_1,\,\dots,\,m_N$. Then elementary algebra yields that both
$P_i(\vec x)$ and $Q_i(\vec x)$ are homogeneous (of the same degree) in the
masses. Since any factor of a homogeneous polynomial is easily seen to be
also homogeneous, we get that the common divisor $Q(\vec x)$ of 
$P_1(\vec x),\,\dots,\,P_s(\vec x)$ is also homogeneous in the variables
$m_1,\dots,m_N$. \qed

\medskip

Our next result, specific to our dynamics, that will be needed later is

\medskip

\subheading{Proposition 3.3} Let $\nu\ge3$, $\harmas$ be a discrete orbit
structure with Property (A) and a $C(N)$-rich symbolic collision sequence
$\Sigma=(\sigma_1,\,\dots,\,\sigma_n)$. Denote by 
$P_1(\vec x),\,\dots,\,P_s(\vec x)$ the polynomials of (4.3) in [S-Sz(1999)]
just as before, and let $Q(\vec x)$ be a common irreducible divisor of 
$P_1(\vec x),\,\dots,\,P_s(\vec x)$. 
(If such a common divisor exists.) Let, finally, 
$\left(\overline{\Sigma},\overline{\Cal A},\vec{\rho}\right)$ be an extended
discrete orbit structure with Property (A) and an extended collision
sequence $\overline{\Sigma}=(\sigma_0,\,\sigma_1,\,\dots,\,\sigma_n)$.
We claim that the irreducible (indecomposable) solution set $V$ of the 
equation $Q(\vec x)=0$ can not even locally coincide with any of the following 
singularity manifolds $C$ defined by one of the following equations:

\medskip

(1) $\Vert v^0_{i_0}-v^0_{j_0}\Vert^2=0$,

\medskip

(2) $\langle v^0_{i_0}-v^0_{j_0};\,\tilde q^0_{i_0}-\tilde q^0_{j_0}
-L\cdot a_0\rangle=0$,

\medskip

(3) $m_{i_0}+m_{j_0}=0$,

\medskip

\noindent
i. e. the irreducible polynomial $Q(\vec x)$ is not equal to any of the
(irreducible) polynomials on the left-hand-sides of (1)---(3). 
(In [S-Sz(1999)] these equations feature Definition 3.18 of the domain
$D\left(\overline{\Sigma},\overline{\Cal A},\vec{\rho}\right)$.) 
Consequently, the open subset 
$V\cap D\left(\overline{\Sigma},\overline{\Cal A},\vec{\rho}\right)$ 
of $V$ is connected and dense in $V$.

\medskip

\subheading{Remark} The first point where we (implicitly) use the condition
$\nu\ge3$ is the irreducibility of the polynomial 
$\Vert v^0_{i_0}-v^0_{j_0}\Vert^2$ on the left-hand-side of (1). Indeed, in the
case $\nu=2$ this polynomial splits as

$$
\aligned
\Vert v^0_{i_0}-v^0_{j_0}\Vert^2\equiv&\left[\left(v^0_{i_0}\right)_1-
\left(v^0_{j_0}\right)_1+\sqrt{-1}\left(\left(v^0_{i_0}\right)_2-
\left(v^0_{j_0}\right)_2\right)\right] \\
\cdot&\left[\left(v^0_{i_0}\right)_1-
\left(v^0_{j_0}\right)_1-\sqrt{-1}\left(\left(v^0_{i_0}\right)_2-
\left(v^0_{j_0}\right)_2\right)\right].
\endaligned
$$
However, it is easy to see that, in the case $\nu\ge3$, all polynomials on the
left of (1)--(3) are indeed irreducible.

\medskip

\subheading{Proof} First of all, we slightly reformulate the negation of the
statement of the proposition. 

Fix one of the three equations of (1)--(3) above, and denote the irreducible
polynomial on its left-hand-side by $R(\vec x)$. By using the quadratic (or
linear) equation $R(\vec x)=0$, we eliminate one variable $x_j$ out of 
$\vec x$ by expressing it as an algebraic function $x_j=g(\vec y)$ of the
remaining variables $\vec y$ of $\vec x$, so that the algebraic function $g$
only contains finitely many field operations and (at most one) square root.
After this elimination $x_j=g(\vec y)$, the meaning of 
$R(\vec x)\equiv Q(\vec x)$ (i. e. the negation of the assertion of the
proposition) is that all algebraic functions
$f_i(\vec x)\equiv\tilde f_i(\vec y)$ in the proof of Lemma 4.2 of 
[S-Sz(1999)] ($i=1,\,\dots,\,s$, constructed for $\harmas$, not for
$\left(\overline{\Sigma},\overline{\Cal A},\vec{\rho}\right)$) are identically
zero in terms of $\vec y$, meaning that every complex orbit segment
$\omega\in\Omega\harmas$, with the initial data $\vec x(\omega)$
in the solution set of $R(\vec x)=0$, is non-sufficient, see also the
``Dichotomy Corollary'' 4.7 in [S-Sz(1999)].
Thus, the negation of the proposition means that
no orbit segment $\omega\in\Omega\harmas$ in the considered singularity
is sufficient.

Now we carry out an induction on the number of balls $N$ quite
in the spirit of the proof of Key Lemma 4.1 of [S-Sz(1999)]. Indeed, the
statement of the proposition is obviously true in the case $N=2$, for in that
case there are no non-sufficient (complex) trajectories 
$\omega\in\Omega\harmas$, i. e. the greatest common divisor of the
polynomials $P_1(\vec x),\,\dots,\,P_s(\vec x)$ is $1$.

Assume now that $N\ge3$, $\nu\ge3$, and the statement of Proposition 3.3
has been proven for all values $N'<N$. Suppose, however, that the statement
is false for some $\harmas$ and extension
$\left(\overline{\Sigma},\overline{\Cal A},\vec{\rho}\right)$
with $N$ balls and Property (A),
i. e. that there exists a common irreducible divisor $Q(\vec x)$ of all the
polynomials $P_1(\vec x),\,\dots,\,P_s(\vec x)$, and $Q(\vec x)$ happens to
be one of the irreducible polynomials on the left-hand-side of (1), (2), or
(3). By using the $C(N)$-richness of $\Sigma=(\sigma_1,\dots,\sigma_n)$,
we select a suitable label $k_0\in\{1,2,\dots,N\}$, say $k_0=N$,
for the substitution 
$m_{N}=0$ along the lines of Lemma 3.1 above, by also ensuring
the existence of the derived schemes $\harmasv$ and
$\left(\overline{\Sigma}',\overline{\Cal A}',\vec{\rho}'\right)$
for the $(N-1)$-ball-system $\{1,2,\dots,N-1\}$
and properties (i)--(iv) (of Lemma 3.1) for $\Sigma'$, see Corollary 4.35
of [S-Sz(1999)] and Lemma 3.1 above. Denote by $\tilde Q(\vec x)$
the polynomial obtained from $Q(\vec x)$ after the substitution $m_{N}=0$.

\medskip

\subheading{Lemma 3.4} The polynomial $\tilde Q(\vec x)$ is not constant.

\medskip

\subheading{Proof} Assume that $\tilde Q(\vec x)\equiv c\in\Bbb C$. The case
$c=0$ means that $m_{N}$ is a divisor of $Q(\vec x)$, thus 
$m_{N}\equiv Q(\vec x)$ (for $Q(\vec x)$ is irreducible), which is 
impossible, since $Q(\vec x)$ has to be one of the polynomials on the 
left-hand-side of (1), (2), or (3). 

If, however, $c\ne0$, then we have that $Q(\vec x)\equiv c+m_{N}S(\vec x)$
with some nonzero polynomial $S(\vec x)$. ($S(\vec x)$ has to be non-zero,
otherwise $Q(\vec x)$ would be a constant, not an irreducible polynomial.)
However, this contradicts to the
proved homogeneity of the polynomial $Q(\vec x)$ in the masses, see 
Proposition 3.2 above. This finishes the proof of the lemma. \qed

\medskip

\subheading{Remark} If one takes a look at the equations (1), (2), (3), he/she
immediately realizes that either $\tilde Q(\vec x)\equiv Q(\vec x)$
(when $Q(\vec x)$ is the polynomial on the left-hand-side of (1) or (2),
or $Q(\vec x)\equiv m_{i_0}+m_{j_0}$ and $N\not\in\{i_0,j_0\}$), or 
$\tilde Q(\vec x)\equiv m_{i_0}$ when $Q(\vec x)\equiv m_{i_0}+m_{j_0}$ and 
$N=j_0$. In this way one can directly and easily verify Lemma 3.4 without
the above ``involved'' algebraic proof. The reason why we still included
the above proof is that later on in this section (in the proof of 
Sub-lemma 3.7) we will need the idea of the presented proof. 

\medskip

The next lemma will use

\subheading{Definition 3.5} Suppose that two
indices $1\le p<q\le n$ and two labels of balls $i,\,j\in\{1,\,\dots,\,N\}$
are given with the additional requirement that if $i=j$, then
$i\in\bigcup_{l=p+1}^{q-1}\sigma_l$. Following the proof of 
Lemma 4.2 of [S-Sz(1999)],
denote by $Q_1(\vec x),\,Q_2(\vec x),\,\dots,\,Q_{\nu}(\vec x)$
($\vec x\in\Bbb C^{(2\nu+1)N+1}$) the polynomials with the property
that for every vector of initial data $\vec x\in D\harmas$ and for every 
$k$, $k=1,\,\dots,\,\nu$, the following equivalence holds true:

$$
\aligned
\left(\exists \, \omega\in\Omega\text{ such that }\vec x(\omega)=
\vec x \; \& \; \left(v_i^p(\omega)\right)_k=\left(v_j^{q-1}(\omega)\right)_k
\right) \\
\Longleftrightarrow Q_k(\vec x)=0.
\endaligned
$$

\medskip

Our next lemma is a strengthened version of Lemma 4.39 of [S-Sz(1999)]:

\medskip

\subheading{Lemma 3.6} Assume that the combinatorial-algebraic scheme
$\harmas$ has Property (A), and use the assumptions and notations of the 
above definition.

We claim that the polynomials 
$Q_1(\vec x),\,Q_2(\vec x),\,\dots,\,Q_{\nu}(\vec x)$ do not have any 
non-constant common divisor. In other words, the equality 
$v_i^p(\omega)=v_j^{q-1}(\omega)$ only takes place on an algebraic variety
with at least two codimensions.

\medskip

\subheading{Remark} Lemma 4.39 of [S-Sz(1999)] asserted that at least one
of the above polynomials $Q_k(\vec x)$ is nonzero. Then, by the permutation
symmetry of the components $k\in\{1,\,\dots,\,\nu\}$, all of these polynomials
are actually nonzero.

\medskip

\subheading{Proof} Induction on the number $N\ge2$.

1. Base of the induction, $N=2$: First of all, by performing the substitution 
$L=0$, we can annihilate all adjustment vectors, see (I), (IV), (VII) of 
Lemma 4.21 in [S-Sz(1999)], and Remark 4.22 there. Then, an elementary 
inspection shows that for any selection of {\it positive real} masses 
$(m_1,\,m_2)$, indeed, the equality $v_i^p(\omega)=v_j^{q-1}(\omega)$ only 
occurs on a manifold with $\nu-1\;(\ge2)$ codimensions in the section
$\Omega\left(\Sigma,\,\Cal A,\,\vec\tau,\,\vec m\right)$ of $\Omega\harmas$ 
corresponding to the selected masses, since any trajectory segment of a
two-particle system with positive masses and $\Cal A=0$ has a very nice,
totally real (and essentially $\nu$-dimensional) representation in the
relative coordinates of the particles: The consecutive, elastic bounces
of a point particle moving uniformly inside a ball of radius $2r$ of
$\Bbb R^\nu$. Therefore, the statement of the lemma is true for $N=2$.

\medskip

Assume now that $N\ge3$, and the lemma has been successfully proven for all 
smaller numbers of balls. By re-labeling the particles, if necessary, we can 
achieve that

\medskip

(i) $N\ne i$, $N\ne j$ and

(ii) if $i=j$, then the ball $i$ has at least one collision between $\sigma_p$
and $\sigma_q$ with a particle different from $N$.

\medskip

For the fixed combinatorial scheme $\harmas$, select a derived scheme 
$\harmasv$ corresponding to the substitution $m_N=0$, see Definition 4.11
and Corollary 4.35 in [S-Sz(1999)].

Our induction step is going to be a proof by
contradiction. Assume, therefore, that the nonzero polynomials
$Q_1(\vec x),\,Q_2(\vec x),\,\dots,\,Q_{\nu}(\vec x)$ do have a common 
irreducible divisor $R(\vec x)$. According to Proposition 3.2 above, the
(irreducible) polynomial $R(\vec x)$ is homogeneous in the variables
$m_1,\,\dots,\,m_N$. Denote by $\tilde R(\vec x)$ the polynomial that we
obtain from $R(\vec x)$ after the substitution $m_N=0$. Similarly to 
Lemma 3.4 above, we claim

\medskip

\subheading{Sub-lemma 3.7} The polynomial $\tilde R(\vec x)$ is not constant.

\medskip

\subheading{Remark} The reason why we cannot simply apply Lemma 3.4 is that
in the proof of that lemma we used the assumption that the irreducible 
polynomial $Q(\vec x)$ was one of the polynomials on the left-hand-side
of (1), (2), or (3) of Proposition 3.3. Right here we do not have such an
assumption.

\medskip

\subheading{Proof} Suppose that $\tilde R(\vec x)\equiv c$, where $c\in\Bbb C$
is a constant, i. e. $R(\vec x)\equiv c+m_N\cdot S(\vec x)$. In the case 
$c=0$ the polynomial $m_N\equiv R(\vec x)$ would be a common divisor of all
the polynomials $Q_1(\vec x),\,Q_2(\vec x),\,\dots,\,Q_{\nu}(\vec x)$, meaning
that in the considered $N$-ball system $\Omega\harmas$ the equation
$m_N(\omega)=0$ implies the equality $v_i^p(\omega)=v_j^{q-1}(\omega)$.
This, in turn, means that in the $(N-1)$-ball system $\{1,\,\dots,\,N-1\}$
(with the discrete algebraic scheme $\harmasv$) the equality
$v_i^p(\omega)=v_j^{q-1}(\omega)$ is an identity, thus contradicting to the
induction hypothesis.

Therefore $c\ne0$, and in the expansion $R(\vec x)\equiv c+m_N\cdot S(\vec x)$
of the irreducible polynomial $R(\vec x)$ we certainly have 
$S(\vec x)\not\equiv 0$, and this means that $R(\vec x)$ is not homogeneous
in the mass variables, thus contradicting to Proposition 3.2. This finishes
the proof of the sub-lemma. \qed

\medskip

\subheading{Finishing the proof of Lemma 3.6} Denote by $\tilde Q_k(\vec x)$
the polynomial that we obtain from $Q_k(\vec x)$ after the substitution
$m_N=0$ ($k=1,\,\dots,\,\nu$), and by $T_k(\vec x)$ the polynomial constructed
for the $(N-1)$-ball system $\{1,\,\dots,\,N-1\}$ (with the discrete algebraic
scheme $\harmasv$) along the lines of Lemma 4.2 of [S-Sz(1999)],
describing the event 
$v_i^p(\omega)=v_j^{q-1}(\omega)$ in this subsystem ($k=1,\,\dots,\,\nu$).
It follows from the induction hypothesis that the zero set $W_k$ of the polynomial
$\tilde Q_k(\vec x)$ (in the phase space of the $(N-1)$-ball system $\harmasv$)
has a codimension-two intersection with the singularities of the
$\harmasv$ system. Indeed, otherwise we would have 
$\left(v_i^p\right)_k\equiv\left(v_j^{q-1}\right)_k$ on some (irreducible)
singularity manifold of the $\harmasv$ subsystem. Then, by the symmetry with
respect to the co-ordinates $k=1,2,\dots,\nu$, we would have $v_i^p\equiv v_j^{q-1}$
on a codimension-one singularity of the subsystem $\harmasv$, contradicting to the
induction hypothesis. This means that the polynomial 
$T_k(\vec x)$ vanishes on the zero set $W_k$ of $\tilde Q_k(\vec x)$,
so the non-constant common divisor $\tilde R(\vec x)$ of 
$\tilde Q_k(\vec x)$ is a common divisor of
$T_1(\vec x),\,\dots,\,T_\nu(\vec x)$, contradicting to the induction
hypothesis. This finishes the proof of Lemma 3.6. \qed

\medskip

\subheading{Continuing the proof of Proposition 3.3}
Denote by $\tilde P_1(\vec x),\,\dots,\, \tilde P_t(\vec x)$ the 
``$P_i$ polynomials'' of the $N$-ball system $\harmas$ with the constraint
$m_{N}=0$ constructed the same way as the polynomials 
$P_1(\vec x),\,\dots\allowmathbreak,P_s(\vec x)$ in (4.3) of [S-Sz(1999)]
for the general case $m_{N}\in\Bbb C$, see also the proof of Lemma 4.2
in the cited paper.
It follows from the algebraic construction of these polynomials that 
the irreducible polynomial $\tilde Q(\vec x)$ is a common divisor of 
$\tilde P_1(\vec x),\,\dots,\, \tilde P_t(\vec x)$.
Recall that, according to our indirect assumption made right before Lemma 3.4,
$Q(\vec x)$ is a common, irreducible divisor of the polynomials
$P_1(\vec x),\dots,P_s(\vec x)$ and, at the same time, $Q(\vec x)$ 
is one of the
polynomials on the left-hand-side of (1), (2), or (3) in Proposition 3.3.
The polynomial $\tilde Q(\vec x)$ was obtained from $Q(\vec x)$ by the 
substitution $m_{N}=0$.

Let us focus now on Lemma 4.9 of [S-Sz(1999)].
The non-sufficiency of the $N$-ball 
system $\harmas$ with the side condition $m_{N}=0$ comes from two sources:
Either from the parallelity of the relative velocities in (2) of Lemma
4.9, or from the non-sufficiency of the $(N-1)$-ball part of the orbit 
segment $\omega$ with the combinatorial scheme $\harmasv$. The first case 
takes place on a complex algebraic set of (at least) $2$ codimensions,
thanks to our original assumption $\nu\ge3$ and Lemma 3.6 above. 
Concerning the
application of the ``non-equality'' Lemma 3.6 above, we note here that once
the velocities $v^p_{i_p}(\omega)$ and $v^{q-1}_{i_q}(\omega)$ are not equal,
the relative velocities $v^p_{N}(\omega)-v^p_{i_p}(\omega)$ and 
$v^{q-1}_{N}(\omega)-v^{q-1}_{i_q}(\omega)=v^p_{N}(\omega)-
v^{q-1}_{i_q}(\omega)$ are not parallel, unless the common velocity
$v^{p}_{N}(\omega)=v^{q-1}_{N}(\omega)$ belongs to the complex line
connecting the different velocities $v^p_{i_p}(\omega)$ and 
$v^{q-1}_{i_q}(\omega)$, which is a codimension-$(\nu-1)$ condition on the
velocity $v^p_{N}(\omega)$. Therefore, the equation $\tilde Q(\vec x)=0$
with the irreducible common divisor $\tilde Q(\vec x)$ of the polynomials 
$\tilde P_1(\vec x),\,\dots,\, \tilde P_t(\vec x)$ can only describe the 
non-sufficiency of the $\harmasv$-part of the system, thus $\tilde Q(\vec x)$ 
should lack the kinetic and mass variables corresponding to the ball with 
label $N$, as the following sub-lemma states:

\medskip

\subheading{Sub-lemma 3.8} The irreducible common divisor 
$\tilde Q(\vec x)$ of the polynomials

$$
\tilde P_1(\vec x),\,\dots,\,\tilde P_t(\vec x)
$$
does not contain the variables with label $N$.

\medskip

\subheading{Proof} Let $D=D\left(\Sigma,\,\Cal A,\,\vec{\tau}\big|m_{N}=0
\right)\subset\Bbb C^{(2\nu+1)N}$ be the open, 
connected and dense domain in $\Bbb C^{(2\nu+1)N}$ defined analogously to
Definition 3.18 (of [S-Sz(1999)]) but incorporating the constraint $m_{N}=0$, 
see also Lemma 3.19 in [S-Sz(1999)].
Let, further, $\Cal N\subset D$ be a small, complex analytic submanifold
of $D$ with complex dimension $(2\nu+1)N-1$, holomorphic to the unit open ball
of $\Bbb C^{(2\nu+1)N-1}$, and such that the polynomial $\tilde Q(\vec x)$ 
($\vec x\in D$) vanishes on $\Cal N$. (Such a manifold $\Cal N\subset D$ 
exists by the induction hypothesis of Proposition 3.3.)
We split the vectors $\vec x\in D$ as $\vec x=(\vec y,\,\vec z)$, so 
that the $\vec z$-part precisely contains the variables bearing the ball label
$N$. We may assume that $\Cal N\subset B_1\times B_2$, where $B_1$ and $B_2$
are small, open balls in the spaces of the components $\vec y$ and $\vec z$,
respectively. 

Assume, to the contrary of the statement of the sub-lemma, that the 
polynomial $\tilde Q(\vec x)\equiv\tilde Q(\vec y,\,\vec z)$ does
depend on the component $\vec z$. Then, for typical but fixed values 
$\vec y_0$ of $\vec y$, the ``slice'' $\left\{\vec y_0\right\}\times B_2$
intersects the manifold $\Cal N$ in a set of complex codimension one. However,
this fact clearly contradicts our earlier observation that the non-sufficiency
of the orbit segments
$\omega\in D=D\left(\Sigma,\,\Cal A,\,\vec{\tau}\big|m_{N}=0\right)$
imposes a codimension-$2$ condition on the coordinates $\vec z$ bearing the 
label $N$. This contradiction finishes the proof of the sub-lemma. \qed

\medskip

\subheading{Finishing the proof of Proposition 3.3}

\medskip

If $Q(\vec x)$ is the left-hand-side of (1) or (2) in 3.3, then we arrive
at the conclusion that the irreducible polynomial 
$\tilde Q(\vec x)\equiv Q(\vec x)$ divides 
$\tilde P_1(\vec x),\,\dots,\, \tilde P_t(\vec x)$, and $N\ne i_0$,
$N\ne j_0$ by Sub-lemma 3.8.
This means, however, that the statement of the proposition is
false for the $(N-1)$-ball system with the discrete algebraic scheme 
$\harmasv$, contradicting to our induction hypothesis.

If, however, the polynomial $Q(\vec x)$ is $m_{i_0}+m_{j_0}$, then in the
case if $N\not\in\{i_0,\,j_0\}$ we arrive at a contradiction just the same
way as above. If $N\in\{i_0,\,j_0\}$, say $N=j_0$, then 
$\tilde Q(\vec x)\equiv m_{i_0}$, and $m_{i_0}$ is a common divisor of
all polynomials $\tilde P_1(\vec x),\,\dots,\, \tilde P_t(\vec x)$ describing
the non-sufficiency of the $\harmasv$ subsystem with the $N-1$ balls 
$\{1,2,\dots,N-1\}$. This means that the above $\harmasv$ 
subsystem is always non-sufficient, provided that $m_{i_0}=0$. In the case
$N\ge4$ it follows from Lemma 4.1 of [S-Sz(1999)]
(applied to the $(N-2)$-ball system
$\{1,2,\dots,N\}\setminus\{i_0,\,N\}$) and from the 
``non-equality'' Lemma 3.6 that almost every $\harmasv$-orbit with $m_{i_0}=0$
is in fact sufficient. One easily checks by inspection
that, in the case $N=3$,  actually every
orbit of the $2$-ball system $\{1,2\}$ with the side condition 
$m_{i_0}=0$ is hyperbolic (sufficient). The obtained contradiction finishes 
the inductive proof of Proposition 3.3. \qed

\medskip

\subheading{Corollary 3.9} Keep all the notations and assumptions of
Proposition 3.3, except that we assume now that 
$\Sigma=(\sigma_1,\,\dots,\,\sigma_n)$ is
$\left(2C(N)+1\right)$-rich and the singularity manifold $C$ is defined by
one of the following equations:

\medskip

(1)' $\Vert v^{k}_{i_k}-v^{k}_{j_k}\Vert^2=0$,

(2)' $\langle v^{k}_{i_k}-v^{k}_{j_k};\,\tilde
q^{k}_{i_k}-\tilde q^{k}_{j_k}-L\cdot a_k\rangle=0$,

(3)' $m_{i_k}+m_{j_k}=0$

\medskip

\noindent
with some $k$, $1\le k\le n$.
Let $Q(\vec x)$ be a common irreducible divisor of the polynomials
$P_1(\vec x),\dots,P_s(\vec x)$ in (4.3) of [S-Sz(1999)]
constructed for the entire discrete structure $\harmas$ as above.

We again claim the same thing: The manifold $C$ and the solution set of
$Q(\vec x)=0$ can not locally coincide. Consequently, an open, dense, and
connected part of the irreducible variety 
$V=\left\{Q(\vec x)=0\right\}$ belongs to the domain
$D\harmas$ of the allowable initial data.

\medskip

\subheading{Proof} We write $\Sigma=(\sigma_1,\dots,\sigma_n)$ in the form
$\Sigma=\left(\Sigma_1,\,\sigma_k,\,\Sigma_2\right)$, where 
$\Sigma_1=(\sigma_1,\dots,\sigma_{k-1})$, 
$\Sigma_2=(\sigma_{k+1},\dots,\sigma_n)$.
Then either $\Sigma_1$ or $\Sigma_2$ is $C(N)$-rich. If $\Sigma_2$ turns out
to be $C(N)$-rich, then we can directly apply the proposition after a simple
time shift $0\longmapsto k$. In the other case, when we only know that
$\Sigma_1$ is $C(N)$-rich, beside the time shift $0\longmapsto k$ an additional
time-reversal is also necessary to facilitate the applicability of the
proposition. \qed

\medskip

Another consequence of Proposition 3.3 is

\subheading{Corollary 3.10} Let $\nu\ge3$, $\harmas$ be a discrete orbit
structure with Property (A) and a $C(N)$-rich symbolic collision sequence
$\Sigma=(\sigma_1,\,\dots,\,\sigma_n)$. Denote by 
$P_1(\vec x),\,\dots,\,P_s(\vec x)$ the polynomials of (4.3) of [S-Sz(1999)]
just as before, and let $Q(\vec x)$ be a common irreducible divisor of 
$P_1(\vec x),\,\dots,\,P_s(\vec x)$. 
(If such a common divisor exists.) Let, finally, 
$\left(\overline{\Sigma},\overline{\Cal A},\vec{\rho}\right)$ be an extended
discrete orbit structure with Property (A) and an extended collision
sequence $\overline{\Sigma}=(\sigma_0,\,\sigma_1,\,\dots,\,\sigma_n)$.
According to Lemma 3.1, we can find a particle, say the one with label $N$, 
and two indices $1\le p<q\le n$ such that

\medskip

(i) $N\in\sigma_p\cap\sigma_q$,

(ii) $N\notin\bigcup_{j=p+1}^{q-1}\sigma_j$,

(iii) $\sigma_p=\sigma_q\Longrightarrow\,\left(\exists j\right)\;
\left(p<j<q\,\&\,\sigma_p\cap\sigma_j\ne\emptyset\right)$, and

(iv) $\Sigma'$ is $(2C(N-1)+1)$-rich on the vertex set $\{1,\dots,N-1\}$.

\noindent
(Here, just as in the case of derived schemes, we denote by $\Sigma'$ the
symbolic sequence that can be obtained from $\Sigma$ by discarding all
edges containing $N$.) Denote by $\tilde Q(\vec x)$ the polynomial that we
obtain from $Q(\vec x)$ after the substitution $m_N=0$. (Obviously,
$\tilde Q(\vec x)\not\equiv0$, otherwise there would not be any sufficient
orbit segment $\omega\in\Omega\left(\Sigma,\,\Cal A,\,\vec\tau\right)$
with $m_N=0$.)

We claim that none of the irreducible polynomials on the left-hand-side of

\medskip

(1) $\Vert v^0_{i_0}-v^0_{j_0}\Vert^2=0$,

\medskip

(2) $\langle v^0_{i_0}-v^0_{j_0};\,\tilde q^0_{i_0}-\tilde q^0_{j_0}
-L\cdot a_0\rangle=0$,

\medskip

(3) $m_{i_0}+m_{j_0}=0$

\medskip
\noindent
is a divisor of $\tilde Q(\vec x)$.

\medskip

\subheading{Proof} Consider and fix an irreducible factor $R(\vec x)$ of
$\tilde Q(\vec x)$. According to Sub-lemma 3.8, the polynomial 
$R(\vec x)\equiv R(\vec y,\,\vec z)$ does not contain any variable bearing 
the label $N$, i. e. $R(\vec x)\equiv R(\vec y)$.

Assume, to the contrary of the statement that we want to prove, that 
the irreducible polynomial $R(\vec x)\equiv R(\vec y)$ is identical to one of
the irreducible polynomials on the left-hand-side of (1), (2), or (3). In 
particular, we have that $N\ne i_0$, $N\ne j_0$. As we saw in the course
of the proof of Sub-lemma 3.8, the algebraic variety 
$V=\left\{\vec y\in\Bbb C^{(2\nu+1)(N-1)+1}\big|\; R(\vec y)=0\right\}$,
defined by one of the equations (1), (2), or (3),
describes the non-sufficiency of the derived system
$\Omega\harmasv$ that one obtains from the original 
$\Omega\harmas$ by taking $m_N=0$, i. e. for any point
$\vec y\in D\harmasv$ it is true that $\vec y\in V$ if and only if there is
some non-sufficient complex orbit segment
$\omega\in\Omega\harmasv$ with $\vec y(\omega)=\vec y$. However, this
statement contradicts to the assertion of Proposition 3.3. \qed

\medskip

\subheading{Remark} Note that the polynomial $\tilde Q(\vec x)$ cannot be
a constant $c\ne0$, otherwise the original polynomial 
$Q(\vec x)=c+m_N\cdot S(\vec x)$ would not be homogeneous in the mass 
variables, see also the proof of Lemma 3.4.

\medskip

The main result of this section is

\medskip

\subheading{Key Lemma 3.11} 
Keep all the notations and notions of this section.
Assume that $\nu\ge3$ and the symbolic collision sequence
$\Sigma=(\sigma_1,\dots,\sigma_n)$ of the discrete algebraic frame
$\harmas$ (with Property (A)) is $C(N)$-rich. 

We claim that all orbit segments $\omega\in\Omega\harmas$ are sufficient
apart from an algebraic variety of codimension-two (at least two, that is),
i. e. the polynomials $P_1(\vec x),\,\dots,\,P_s(\vec x)$ of (4.3) of
[S-Sz(1999)] do not have a non-constant common divisor.

\subheading{Proof} The inductive proof employs many of the ideas of the proof
of Proposition 3.3 and it will use the statement of the proposition itself.
(More precisely, the statement of Corollary 3.9 is to be used.)

Indeed, the assertion of this key lemma is trivially true in the case $N=2$,
for in that case there are no non-sufficient, complex orbit segments
$\omega\in\Omega\harmas$ at all.

Assume that $N\ge3$, and the statement of the key lemma has been successfully
proven for all smaller values ($2\le$) $N'<N$. Our induction step is going
to be a proof by contradiction. Suppose, therefore, that the polynomials
$P_1(\vec x),\,\dots,\,P_s(\vec x)$ do have a common irreducible divisor
$Q(\vec x)$. Following the assertion of Lemma 3.1, select a suitable label
$k_0\in\{1,\,\dots,\,N\}$ for the substitution $m_{k_0}=0$ so that a derived
scheme $\harmasv$ (with Property (A)) exists for the arising $(N-1)$-ball
system $\{1,\,\dots,\,N\}\setminus\{k_0\}$ with a symbolic sequence
$\Sigma'$, possessing the properties (1)--(4) of Lemma 3.1, the same way as
we did in the proof of Proposition 3.3. Without loss of generality,
we may assume that $k_0=N$.

Consider now the original system $\Omega\harmas$ with the constraint
$m_{N}=0$. After the substitution $m_{N}=0$ the polynomial $Q(\vec x)$
becomes a new, non-constant polynomial $\tilde Q(\vec x)$, see the proof of
Sub-lemma 3.7 above. Let $S(\vec x)$ be an irreducible divisor of 
$\tilde Q(\vec x)$. The (indecomposable) algebraic variety
$V=\left\{S(\vec x)=0\right\}$ has one codimension in the submanifold 
$\tilde\Omega=\Omega_{m_N=0}$ of $\Omega\harmas$, and for every $\vec x\in V$
there exists a non-hyperbolic complex orbit segment $\omega\in\Omega\harmas$
with $\vec x(\omega)=\vec x$ and $m_{N}(\omega)=0$.
As far as the non-sufficiency of
the orbits $\omega\in\Omega\harmas$ with $m_{N}(\omega)=0$ is
concerned, we again take a close look at Lemma 4.9 of [S-Sz(1999)].
We saw earlier (see the proof of Lemma 3.6, which clearly carries over to the
models subjected to the side condition $m_N=0$) that the parallelity of the 
relative velocities $v^p_{N}(\omega)-v^p_{i_p}(\omega)$ and
$v^{q-1}_{N}(\omega)-v^{q-1}_{i_q}(\omega)$ takes place on a manifold with
codimension at least two in our case of $\nu\ge3$. Therefore, according to 
Lemma 4.9 of [S-Sz(1999)], the ``codimension-one event'' 
$\vec x(\omega)\in V$ ($\Longleftrightarrow S(\vec x(\omega))=0$),
for orbits with $m_{N}(\omega)=0$, can only be equivalent to the 
non-sufficiency of the $\{1,\,\dots,\,N-1\}$-part
$\text{trunc}(\omega)\in\Omega\harmasv$ of the system. In this way
it follows from Lemma 4.9 that the irreducible polynomial $S(\vec x)$
lacks all variables bearing the label $N$, see also the statement and the
proof of Sub-lemma 3.8. We conclude that for every $\vec x\in V$ (i. e. with
$S(\vec x)=0$) there exists an orbit segment $\omega\in\Omega\harmas$ with
$m_{N}(\omega)=0$, $\vec x(\omega)=\vec x$, and a non-sufficient truncated 
segment $\omega'=\text{trunc}(\omega)\in\Omega\harmasv$.
According to Proposition 3.3 above (applied to the 
$(N-1)$-ball system $\{1,\,\dots,\,N-1\}$ with the algebraic
scheme $\harmasv$), the variety $\left\{S(\vec x)=0\right\}$ does not even
locally coincide with the singularities of the complex dynamics
$\Omega\harmasv$. This means that a codimension-one family of complex orbit 
segments $\omega'=\text{trunc}(\omega)\in\Omega\harmasv$, 
$\vec x(\omega)\in V$, is not sufficient. This, in turn, contradicts the 
induction hypothesis of the proof of Key Lemma 3.11 by actually finishing it.
\qed

\bigskip \bigskip

\heading
4. Finishing the Proof of Ergodicity
\endheading

\bigskip \bigskip

\subheading{From $\Bbb C$ back to $\Bbb R$} 

\medskip

First of all, we transfer the main
result of the previous section (Key Lemma 3.11) from the complex set-up back
to the real case. This result will be an almost immediate consequence of
Key Lemma 3.11.

Fix a discrete algebraic scheme $\harmas$ for $N$ balls with Property (A)
(see Definition 3.31 in [S-Sz(1999)])
and a $C(N)$-rich symbolic collision sequence
$\Sigma=(\sigma_1,\dots,\sigma_n)$. (The definition of the threshold $C(N)$
can be found in Lemma 3.1.) Denote by
$\Omega_{\Bbb R}=\Omega_{\Bbb R}\harmas$ the set of all elements
$\omega\in\Omega\harmas$ for which

\medskip

(1) all kinetic functions $\left(\tilde q_i^k(\omega)\right)_j$,
$\left(v_i^k(\omega)\right)_j$, $m_i(\omega)$, and $L(\omega)$ take real
values, $i=1,\dots,N$; $k=0,1,\dots,n$; $j=1,\dots,\nu$;

(2) $\tau_k(\omega)=t_k(\omega)-t_{k-1}(\omega)>0$ for $k=1,\dots,n$;

(3) out of the two real roots of (3.8) of [S-Sz(1999)]
the root $\tau_k$ is always selectedas the smaller one, $k=1,\dots,n$.

\medskip

It is clear that either $\Omega_{\Bbb R}=\Omega_{\Bbb R}\harmas$ is a
$\left((2\nu+1)N+1\right)$-dimensional, real analytic submanifold of
$\Omega=\Omega\harmas$, or $\Omega_{\Bbb R}=\emptyset$. Of course, we will
never investigate the case $\Omega_{\Bbb R}=\emptyset$.

Consider the corresponding polynomials $P_1(\vec x),\dots,P_s(\vec x)$ of
(4.3) of [S-Sz(1999)] describing the non-sufficiency of the complex orbit 
segments $\omega\in\Omega\harmas$, along the lines of Lemma 4.2 of 
[S-Sz(1999)], in terms of the kinetic data 
$\vec x=\vec x(\omega)$. According to the statement in the third paragraph
on p. 61 of [S-Sz(1999)], these polynomials
$P_i(\vec x)$ admit real coefficients. By Key Lemma 3.11, the greatest common
divisor of $P_1(\vec x),\dots,P_s(\vec x)$ is $1$, hence the common zero set

$$
\left\{\vec x\in\Bbb R^{(2\nu+1)N+1}\big|\; P_1(\vec x)=P_2(\vec x)=\dots
=P_s(\vec x)=0\right\}
$$
of these polynomials does not contain any smooth real submanifold of
(real) dimension $(2\nu+1)N$. In this way we obtained

\medskip

\subheading{Proposition 4.1} Use all the notions, notations and assumptions
from above. There exists no smooth, real submanifold $\Cal M$ of 
$\Omega_{\Bbb R}$ with 
$\text{dim}_{\Bbb R}\Cal M=\text{dim}_{\Bbb R}\Omega_{\Bbb R}-1
\left(=(2\nu+1)N\right)$ and with the property that all orbit segments
$\omega\in\Cal M$ are non-sufficient. (For the concept of non-sufficiency, 
please see \S2.) \qed

\bigskip

\subheading{The ``Fubini-type'' Argument} 

\medskip

Our dynamics $\Omega\harmas$ has the obvious feature that the variables
$m_i=m_i(\omega)$ ($i=1,\dots,N$) and $L(\omega)$ (the so called outer
geometric parameters) remain unchanged during the time-evolution. Quite
naturally, we do not need Proposition 4.1 directly but, rather, we need to
use its analog for (almost) every fixed $(N+1)$-tuple
$(m_1,\dots,m_N;L)\in\Bbb R^{N+1}$. This will be easily achieved by a
classical ``Fubini-type'' product argument. The result is

\medskip

\subheading{Proposition 4.2} Use all the notions, notations and assumptions
from above. Denote by

$$
NS=NS\harmas=\left\{\omega\in\Omega_{\Bbb R}\harmas\big|\;
\text{dim}_{\Bbb C}\Cal N(\omega)>\nu+1\right\}
$$
the set of all non-sufficient, real orbit segments 
$\omega\in\Omega_{\Bbb R}=\Omega_{\Bbb R}\harmas$. (For the definition of the
complex neutral space $\Cal N(\omega)$, please see (3.21) in [S-Sz(1999)].) 
Finally, we use the notation

$$
\Omega_{\Bbb R}(\vec m,\,L)=\left\{\omega\in\Omega_{\Bbb R}\big|\;
\vec m(\omega)=\vec m, \text{ and } L(\omega)=L\right\}
$$
for any given $(N+1)$-tuple
$(\vec m,\,L)=(m_1,\dots,m_N,\,L)\in\Bbb R^{N+1}$. We claim that for almost
every $(\vec m,\,L)\in\Bbb R^{N+1}$ (for which 
$\Omega_{\Bbb R}\left(\vec m,\,L\right)\ne\emptyset$) the
intersection $NS\cap\Omega_{\Bbb R}(\vec m,\,L)$ has at least $2$ 
codimensions in $\Omega_{\Bbb R}(\vec m,\,L)$.

\medskip

\subheading{Remark 4.3} As it is always the case with such algebraic systems,
the exceptional zero-measure set of the parameters $(\vec m,\,L)$ turns out
to be a countable union of smooth, proper submanifolds of $\Bbb R^{N+1}$.

\medskip

\subheading{Proof of Proposition 4.2}
It is clear that the statement of the proposition
is a local one, therefore it is enough to prove that for any small, open 
subset $U_0\subset\Omega_{\Bbb R}$ of 
$\Omega_{\Bbb R}=\Omega_{\Bbb R}\harmas$ the set

$$
\left\{(\vec m,\,L)\in\Bbb R^{N+1}\big|\; \text{dim}_{\Bbb R}\left(
NS\cap\Omega_{\Bbb R}(\vec m,\,L)\cap U_0\right)\ge 2\nu N-1\right\}
$$
of the ``bad points'' $(\vec m,\,L)$ has zero Lebesgue measure. The points
$\omega\in U_0$ can be identified locally (in $U_0$) with the vector 
$\vec x=\vec x(\omega)\in D_{\Bbb R}=D\harmas\cap\Omega_{\Bbb R}$ of their
initial coordinates. After this identification the small open set 
$U_0\subset\Omega_{\Bbb R}$ naturally becomes an open subset 
$U_0\subset D_{\Bbb R}$. Furthermore, we split the points $\vec x\in U_0$
as $\vec x=\left((\vec m,\,L),\,\vec y\right)$, where $\vec y$ contains all
variables other than $m_1,\,\dots,\,m_N,\,L$. In this way we may assume that
$U_0$ has a product structure $U_0=B_0\times B_1$ of two small open balls,
so that $B_0\subset\Bbb R^{N+1}$, while the open ball 
$B_1\subset\Bbb R^{2\nu N}$ contains the $\vec y$-parts of the points
$\vec x=\left((\vec m,\,L),\,\vec y\right)\in U_0$. 

Assume that the statement of the proposition is false. Then there exists
a small open set $U_0=B_0\times B_1\subset\Bbb R^{N+1}\times\Bbb R^{2\nu N}$
(with the above splitting) and there is a positive number $\epsilon_0$ such
that the set

$$
\aligned
A_0=&\big\{(\vec m,\,L)\in B_0\big|\; \left((\vec m,\,L)\times B_1\right)
\cap NS\text{ contains a} \\
&(2\nu N-1)\text{-dimensional, smooth, real submanifold with inner radius }
>\epsilon_0\big\}
\endaligned
$$
has a positive Lebesgue measure in $B_0$. Then one can find an orthogonal
projection $P:\; \Bbb R^{2\nu N}\to H$ onto a hyperplane $H$ of
$\Bbb R^{2\nu N}$ such that, by taking
$\Pi(\vec x)=\Pi\left((\vec m,\,L),\,\vec y\right)=P(\vec y)$,
($\Pi:\; \Bbb R^{(2\nu+1)N+1}\to H$), the set 

$$
\aligned
A_1=&\big\{(\vec m,\,L)\in B_0\big|\; \Pi\left[\left((\vec m,\,L)
\times B_1\right)\cap NS\right]\text{ contains} \\
&\text{an open ball of radius }>\epsilon_0/2\text{ in } H\big\}
\endaligned
$$
has positive Lebesgue measure in $B_0$. By the Fubini Theorem the set

$$
\tilde\Pi\left[NS\cap(B_0\times B_1)\right]
$$
has positive Lebesgue measure in $B_0\times H$, where

$$
\tilde\Pi(\vec x)=\tilde\Pi\left((\vec m,\,L),\,\vec y\right)=
\left((\vec m,\,L),\,P(\vec y)\right)\in B_0\times H
$$
for $\vec x\in B_0\times B_1$. However,
$\text{dim}_{\Bbb R}(B_0\times H)=(2\nu+1)N$, and
$\text{dim}_{\Bbb R}\left(NS\cap\Omega_{\Bbb R}\right)\le(2\nu+1)N-1$
(according to Proposition 4.1). Thus, we obtained that the real algebraic set

$$
\tilde\Pi\left(NS\cap(B_0\times B_1)\right)\subset B_0\times H
$$
has dimension strictly less than 
$\text{dim}_{\Bbb R}(B_0\times H)=(2\nu+1)N$, yet it has positive Lebesgue 
measure in the space $B_0\times H$. The obtained contradiction finishes the
proof of Proposition 4.2. \qed

\bigskip

\heading
Finishing the Proof of the Theorem
\endheading

\bigskip

We will carry out an induction with respect to the number of balls $N$ 
($\ge2$). For $N=2$ the system is well known to be a strictly dispersive
billiard flow (after the obvious reductions $m_1v_1+m_2v_2=0$,
$m_1||v_1||^2+m_2||v_2||^2=1$ ($m_1,\,m_2>0$), and after the factorization
with respect to the uniform spatial translations, as usual) and, as such,
it is proved to be ergodic by Sinai in [Sin(1970)], see also the paper 
[S-W(1989)] about the case of different masses.

Assume now that $N\ge3$, $\nu\ge3$, and the theorem has been successfully 
proven for all smaller numbers of balls $N'<N$. Suppose that a billiard flow

$$
\flow=\left(\bold M_{\vec m,L},\,\left\{S^t_{\vec m,L}\right\},\,
\mu_{\vec m,L}\right)
$$
is given for $N$ balls and outer geometric parameters
$(\vec m,\,L)=(m_1,\,\dots,\,m_N,\,L)$ ($m_i>0$, $L>0$) in such a way that,
besides the always assumed properties (2.1.1)--(2.1.2),

\medskip

(*) the vector $(\vec m,\,L)$ of geometric parameters is such that for any
subsystem $1\le i_1<i_2<\dots<i_{N'}\le N$ ($2\le N'\le N$) and for any
$C(N')$-rich discrete algebraic scheme $\harmas$ (with Property (A)), for this
subsystem $\left(m_{i_1},\,\dots,\,m_{i_{N'}},\,L\right)$ it is true that the
parameter vector $\left(m_{i_1},\,\dots,\,m_{i_{N'}},\,L\right)$ does not 
belong to the zero-measured exceptional set of parameters featuring
Proposition 4.2.

\medskip

According to Lemma 4.1 of [K-S-Sz(1990)], the set $R_2\subset\bold M$ of the
phase points with at least two singularities on their trajectories is a
countable union of smooth submanifolds of $\bold M$ with codimension two,
so this set $R_2$ can be safely discarded in the proof, for it is slim, see
also \S2 about the slim sets. Secondly, by the induction hypothesis and by
Theorem 5.1 of [Sim(1992-I)] (adapted to the case of different masses) there
is a slim subset $S_1\subset\bold M$ such that for every phase point
$x\in\bold M\setminus S_1$ 

\medskip

(i) $S^{(-\infty,\infty)}x$ contains at most one singularity, and

\medskip

(ii) $S^{(-\infty,\infty)}x$ contains an arbitrarily large number of 
consecutive, connected collision graphs.

\medskip

\noindent
(In the case of a singular trajectory $S^{(-\infty,\infty)}x$ we require that
both branches contain an arbitrarily large number of consecutive, connected 
collision graphs.) Then, by Proposition 4.2 just proved, there is another
slim subset $S_2\supset S_1$ of $\bold M$ such that

\medskip

(H) for every $x\in\bold M\setminus S_2$ the trajectory $S^{(-\infty,\infty)}x$
contains at most one singularity and it is sufficient (or, geometrically
hyperbolic).

\medskip

According to Theorem 6.1 of [Sim(1992-I)] (easily adapted to the case of 
different masses) and Proposition 4.2, the so called Chernov-Sinai Ansatz
(see \S2) holds true, i. e. for almost every singular phase point 
$x\in\Cal S\Cal R^+$ the positive semi-trajectory $S^{(0,\infty)}x$ is
non-singular and sufficient.

This is the point where the Fundamental Theorem for Algebraic Semi-dispersive
Billiards (Theorem 4.4 in [B-Ch-Sz-T(2002)]) comes to play! According to that
theorem, by also using the crucial conditions (H) and the Ansatz above, it is
true that for every phase point
$x\in\left(\text{int}\bold M\right)\setminus S_2$ some open neighborhood 
$U_x$ of $x$ in $\bold M$ belongs to a single ergodic component of the 
considered billiard flow $\left\{S_{\vec m,L}^t\right\}$. Since the set
$\left(\text{int}\bold M\right)\setminus S_2$ contains an arc-wise connected
set $C$ with full $\mu$-measure (see Proposition 2.7.9 above),
we get that the entire set $C$ belongs to a single ergodic component
of the flow $\left\{S_{\vec m,L}^t\right\}$. This finishes the proof of the
ergodicity theorem. \qed

\bigskip \bigskip

\subheading{Concluding Remark: The Irrational Mass Ratio}

\bigskip

Due to the natural reduction $\sum_{i=1}^N m_iv_i=0$ (which we always assume),
in \S\S1--2 we had to factorize the configuration space with respect to
spatial translations: $(q_1,\dots,q_N)\sim(q_1+a,\dots,q_N+a)$ for all
$a\in\Bbb T^\nu$. It is a remarkable fact, however, that (despite the 
reduction $\sum_{i=1}^N m_iv_i=0$) even without this translation factorization
the system still retains the Bernoulli mixing property, provided that the
masses $m_1,\dots,m_N$ are rationally independent. (We note that dropping the
above mentioned configuration factorization obviously introduces $\nu$ zero
Lyapunov exponents.) For the case $N=2$ (i. e. two disks)
this was proven in [S-W(1989)] by successfully applying D. Rudolph's following
theorem on the B-property of isometric group extensions of Bernoulli shifts 
[R(1978)]:

Suppose that we are given a dynamical system $(M,T,\mu)$ with a probability
measure $\mu$ and an automorphism $T$. Assume that a compact metric group
$G$ is also given with the normalized Haar measure $\lambda$ and left invariant
metric $\rho$. Finally, let $\varphi\colon\; M\to G$ be a measurable map.
Consider the skew product dynamical system $(M\times G,S,\mu\times\lambda)$ 
with $S(x,g)=\left(Tx,\varphi(x)\cdot g\right)$, 
$x\in M$, $g\in G$. We call the system $(M\times G,S,\mu\times\lambda)$ an 
isometric group extension of the base (or factor) $(M,T,\mu)$. (The phrase
``isometric'' comes from the fact that the left translations 
$\varphi(x)\cdot g$ are isometries of the group $G$.) Rudolph's mentioned
theorem claims that the isometric group extension 
$(M\times G,S,\mu\times\lambda)$ enjoys the B-mixing property as long as it is
at least weakly mixing and the factor system $(M,T,\mu)$ is a B-mixing system.

But how do we apply this theorem to show that the typical
system of $N$ hard balls
in $\Bbb T^\nu$ with $\sum_{i=1}^N m_iv_i=0$ is a Bernoulli flow, even if we
do not make the factorization (of the configuration space) with respect to
uniform spatial translations? It is simple. The base system $(M,T,\mu)$
of the isometric group extension $(M\times G,S,\mu\times\lambda)$ will be the
time-one map of the factorized (with respect to spatial translations) hard 
ball system. The group $G$ will be just the container torus $\Bbb T^\nu$
with its standard Euclidean metric $\rho$ and normalized Haar measure
$\lambda$. The second component $g$ of a phase point
$y=(x,g)\in M\times G$ will be just the position of the center of the (say) 
first ball in $\Bbb T^\nu$. Finally, the governing translation
$\varphi(x)\in\Bbb T^\nu$ is quite naturally the total displacement

$$
\int\Sb 0\endSb\Sp 1\endSp v_1(x_t)dt \qquad (\text{mod }\Bbb Z^\nu)
$$
of the first particle while unity of time elapses. In the previous sections
the B-mixing property of the factor map $(M,T,\mu)$ has been proven
successfully for typical geometric parameters $(m_1,\dots,m_N;\, L)$.
Then the key step in proving the B-property of the isometric
group extension $(M\times G,S,\mu\times\lambda)$ is to show that the latter
system is weakly mixing. This is just the essential contents of the
article [S-W(1989)], and it takes advantage of the assumption of rational
independence of the masses. Here we are only presenting to the reader
the outline of that proof. As a matter of fact, we not only proved
the weak mixing property of the extension $(M\times G,S,\mu\times\lambda)$,
but we showed that this system has in fact the K-mixing property by proving
that the Pinsker partition $\pi$ of $(M\times G,S,\mu\times\lambda)$ is
trivial. (The Pinsker partition is, by definition, the finest invariant,
measurable partition of the dynamical 
system with respect to which the factor system
has zero metric entropy. A dynamical system is K-mixing if and only if
its Pinsker partition is trivial, i. e. it consists of only the sets with
measure zero and one, see [K-S-F(1980)].) In order to show that the Pinsker
partition is trivial, in [S-W(1989)] we constructed a pair of measurable
partitions $(\xi^s,\,\xi^u)$ for $(M\times G,S,\mu\times\lambda)$ made up
by open and connected sub-manifolds of the local stable and unstable 
manifolds, respectively. It followed by standard methods (see [Sin(1968)]) 
that the partition $\pi$ is coarser than each of $\xi^s$ and $\xi^u$. Due to 
the $S$-invariance of $\pi$, we have that $\pi$ is coarser than

$$
\bigwedge\Sb n\in\Bbb Z\endSb S^n\xi^s\wedge
\bigwedge\Sb n\in\Bbb Z\endSb S^n\xi^u.
\tag *
$$
In the final step, by using now the rational independence of the masses,
we showed that the partition in $(*)$ is, indeed, trivial.

\bye